\documentclass{article}

\usepackage{amssymb,amsmath,amscd,dsfont}
\newcommand{\op}[1]{\ensuremath{\operatorname{#1}}}
\newcommand{\ol}[1]{\ensuremath{\overline{#1}}}
\newcommand{\ul}[1]{\ensuremath{\underline{#1}}}
\newcommand{\fa}{\ensuremath{\;\text{\textup{ for all }}\;}}
\newcommand{\cC}{\ensuremath{\mathcal{C}}}
\newcommand{\cD}{\ensuremath{\mathcal{D}}}
\newcommand{\cG}{\ensuremath{\mathcal{G}}}
\newcommand{\cH}{\ensuremath{\mathcal{H}}}
\newcommand{\cK}{\ensuremath{\mathcal{K}}}
\newcommand{\cM}{\ensuremath{\mathcal{M}}}
\newcommand{\cP}{\ensuremath{\mathcal{P}}}
\newcommand{\cQ}{\ensuremath{\mathcal{Q}}}
\newcommand{\cV}{\ensuremath{\mathcal{V}}}
\newcommand{\fK}{\ensuremath{\mathfrak{K}}}
\newcommand{\fg}{\ensuremath{\mathfrak{g}}}
\newcommand{\fG}{\ensuremath{\mathfrak{G}}}
\newcommand{\fh}{\ensuremath{\mathfrak{h}}}
\newcommand{\fl}{\ensuremath{\mathfrak{l}}}
\newcommand{\mf}[1]{\ensuremath{\mathfrak{#1}}}
\newcommand{\R}{\ensuremath{\mathds{R}}}
\newcommand{\N}{\ensuremath{\mathds{N}}}
\newcommand{\Z}{\ensuremath{\mathds{Z}}}
\newcommand{\id}{\ensuremath{\operatorname{id}}}
\newcommand{\pr}{\ensuremath{\operatorname{pr}}}
\newcommand{\ev}{\ensuremath{\operatorname{ev}}}
\newcommand{\Ad}{\ensuremath{\operatorname{Ad}}}
\newcommand{\Obj}{\ensuremath{\operatorname{Obj}}}
\newcommand{\Mor}{\ensuremath{\operatorname{Mor}}}
\newcommand{\Out}{\ensuremath{\operatorname{Out}}}
\newcommand{\cAd}{\ensuremath{\operatorname{\mathcal{A}\mathnormal{d}}}}
\newcommand{\cExp}{\ensuremath{\operatorname{\mathcal{E}\mathnormal{xp}}}}
\newcommand{\cMor}{\ensuremath{\operatorname{\mathcal{M}\mathnormal{or}}}}
\newcommand{\cFun}{\ensuremath{\operatorname{\mathcal{F}\mathnormal{un}}}}
\newcommand{\cBun}{\ensuremath{\operatorname{\mathcal{B}\mathnormal{un}}}}
\newcommand{\cAut}{\ensuremath{\operatorname{\mathcal{A}\mathnormal{ut}}}}
\newcommand{\Aut}{\ensuremath{\operatorname{Aut}}}
\newcommand{\der}{\ensuremath{\operatorname{der}}}
\newcommand{\Diff}{\ensuremath{\operatorname{Diff}}}
\newcommand{\Gau}{\ensuremath{\operatorname{Gau}}}
\newcommand{\Hom}{\ensuremath{\operatorname{Hom}}}
\newcommand{\tx}[1]{\ensuremath{\text{#1}}}
\newcommand{\se}{\ensuremath{\nobreak\subseteq\nobreak}}
\newcommand{\from}{\ensuremath{\nobreak:\nobreak}}
\renewcommand{\to}{\ensuremath{\nobreak\rightarrow\nobreak}}
\usepackage{amssymb,amsmath,amscd}

\usepackage[amsmath,thmmarks]{ntheorem}
\theoremseparator{.}
\theoremsymbol{\rule{1ex}{1ex}}
\theorembodyfont{\upshape}
\newtheorem{definition}{Definition}[section]
\newtheorem{remark}[definition]{Remark}

\newtheorem{example}[definition]{Example}

\newtheorem*{proof}{Proof}

\theorembodyfont{\itshape}

\newtheorem*{nntheorem}{Theorem}
\theoremsymbol{}
\newtheorem{lemma}[definition]{Lemma}
\newtheorem{proposition}[definition]{Proposition}
\newtheorem{theorem}[definition]{Theorem}
\newtheorem{corollary}[definition]{Corollary}

\usepackage[notref,notcite]{showkeys}
\usepackage[plainpages=false,pdfpagelabels]{hyperref}

\renewcommand{\Obj}{\ensuremath{\op{Ob}}}
\renewcommand{\fG}{\ensuremath{\mathfrak{G}}}
\newcommand{\idObject}{\ensuremath{\mathds{1}}}

\usepackage[all,cmtip,line]{xy}
\CompileMatrices
\newcommand{\blabel}[1]{\setlength\fboxrule{.0em}\fbox{\ensuremath{\scriptstyle #1}}}

\begin{document}
\sloppy 

\title{Principal 2-bundles and their gauge 2-groups}
\author{Christoph Wockel\footnote{ Mathematisches Institut,
Bunsenstraße 3-5, 37073 G\"ottingen, Germany}\\[.5em]\texttt{\small
christoph@wockel.eu}} \date{} \maketitle

\begin{abstract}
In this paper we introduce principal 2-bundles
and show how they are classified by non-abelian \v{C}ech
cohomology. Moreover, we show that their gauge 2-groups
can be described by
2-group-valued functors, much like in classical bundle theory. Using this,
we show that, under some mild requirements, these gauge 2-groups possess a
natural smooth structure. In the last section we provide some explicit
examples.\\[\baselineskip] \textbf{MSC:} 55R65, 22E65, 81T13
\end{abstract}

\section*{Introduction}

This paper gives a precise description of globally defined geometric
objects, which are classified by non-abelian \v{C}ech cohomology. The
general philosophy is to realise these geometric structures as
categorified principal bundles, i.e., principal bundles, where sets are
replaced by categories, maps by functors and commuting diagrams by
natural equivalences of functors (satisfying canonical coherence
conditions). The control on these categorified geometric structures is
the amount of natural transformations that one allows to differ from
the identity. For instance, allowing only identity transformations in
the definition of a categorified group below makes it accessible as a
crossed module and thus in terms of ordinary group theory.

The main enrichment that comes from categorification is the existence
of ``higher morphisms'' between morphisms, that are not present in
the set-theoretical setup. Two prominent examples amongst are
homotopies between continuous maps and bimodule morphisms between
bimodules (as morphisms between rings or $C^{*}$-algebras). These
higher morphisms lead to very rich structure, because
it allows a more flexible concept of invertible morphisms, namely
invertibility ``up to higher morphisms''. In the two examples
mentioned, these are the concepts of homotopy equivalence and Morita
equivalence.

The main idea of this paper is to represent non-abelian cohomology
classes as semi-strict principal 2-bundles, i.e., smooth 2-spaces with
a locally trivial strict action of a strict Lie 2-group. This
particular subclass of principal 2-bundles is quite easily accessible,
while the theory in its full generality is much more involved (cf.\
\cite{bartels06HigherGaugeTheoryI:2-Bundles}). In the first section, we
develop the concept of principal 2-bundles from first principles and
show very precisely in the second section, how semi-strict principal
2-bundles are classified by non-abelian \v{C}ech cohomology:

\begin{nntheorem}
 If $\cG$ is a strict Lie 2-group and $M$ is a smooth manifold, then
 semi-strict principal $\cG$-2-bundles over $M$ are classified up to
 Morita equivalence by $\check{H}(M,\cG)$.
\end{nntheorem}

In addition, we provide a geometric way to think of the band of a
semi-strict principal 2-bundle. The initial ideas exposed in this
sections are not new, the earliest reference we found was
\cite{dedecker60SurLaCohomologieNonAbelienne}. What is in scope is a
clear and down to earth development of the idea in order to
make the subject easily accessible.

The paper also aims at opening the subject to infinite-dimensional
Lie theory, and we completely neglect the gauge theoretic motivation of
the theory (cf.\ \cite{baezSchreiber05HigherGaugeTheory},
\cite{schreiberWaldorf07ParallelTransportAndFunctors} and
\cite{SchreiberWaldorf08Connections-on-non-abelian-Gerbes-and-their-H%
olonomy} for this). The treatment of symmetry groups of principal
2-bundles in the third section is a first example. There we show how to
identify the automorphism 2-group of a semi-strict principal 2-bundle
with 2-group valued functors and use this to put smooth
structures on these automorphism 2-groups. The main result of
the third section is the following. 

\begin{nntheorem}
 Let $\cP\to M$ be a semi-strict principal $\cG$-2-bundle. Assume that
 $M$ is compact, that $\cG$ is locally exponential, and that the action
 of $\cG$ on $\cP$ is principal. Then $\cC^{\infty}(\cP,\cG_{\Ad})^{\cG}$ is a
 locally exponential strict Lie 2-group with strict Lie 2-algebra
 $\cC^{\infty}(\cP,L(\cG)_{\op{ad}})^{\cG}$.
\end{nntheorem}

In the classical setup of principal bundles, given as smooth manifolds
on which Lie groups act locally trivially, this isomorphism of gauge
transformations and group valued functions was the dawning of the
global formulation of gauge theories in terms of principal bundles,
which lead to its fancy developments.

The last section treats examples, in particular bundle gerbes (or
groupoid extensions, much like \cite[Sect.\ 
2]{Laurent-GengouxStienonXu09Non-abelian-differentiable-gerbes},
\cite[Sect.\
3]{ginotStienon08G-gerbesPrincipal2-bundlesAndCharacteristicClasses},
\cite[Sect.\
4]{moerdijk02IntroductionToTheLanguageOfStacksAndGerbes}). This
section is not exhaustive, it should give an intuitive idea for
relating bundle gerbes and principal 2-bundles and give some further
examples. In the end, we provided a short appendix with some basic
concepts of locally convex Lie theory.

In comparison to many other expositions of this subject
the principal bundles that we consider are globally defined objects, 
considered as smooth 2-spaces, together with a locally trivial smooth action of a 2-group.
Most of the other approaches consider more general base spaces than we
do by allowing hypercovers, arbitrary surjective submersions or even
more general Lie groupoids. This text only treats \v{C}ech groupoids as
surjective submersions. Some global aspects of constructing higher
bundles can be found in \cite{Laurent-GengouxStienonXu09Non-abelian-differentiable-gerbes}, \cite[Ex.\
2.19]{ginotStienon08G-gerbesPrincipal2-bundlesAndCharacteristicClasses},
\cite[Th.\ 3.1]{moerdijk02IntroductionToTheLanguageOfStacksAndGerbes}
and \cite[Sect.\ 2.7]{breen92OnTheClassificationOf2-GerbesAnd2-Stacks},
part of which served as a motivation for our construction. However,
there is no reference known to the author, which treats principal
2-bundles as smooth \textit{2-spaces with a locally trivial action} of
a structure 2-group. The objects that come closest to what we
have in mind are bundle gerbes, as treated in
\cite{murray96BundleGerbes},
\cite{aschieriCantiniLuigiBranislav05NonabelianBundleGerbesTheirDifferentialGeomoetryAndGaugeTheory},
\cite{schweigertWaldorf07GerbesAndLieGroups}. Another global way for
describing principal bundles is in terms of their transport 2-functors
(cf.\ \cite{schreiberWaldorf07ParallelTransportAndFunctors} for
ordinary principal bundles and \cite{SchreiberWaldorf08Connections-on%
-non-abelian-Gerbes-and-their-Holonomy} for principal 2-bundles), which
is closely related to our approach.

Although there frequently exist more systematic approaches to the
things we present, we avoid the introduction of more general
frameworks (such as internal categories, 2-categories, etc.). This is
done to keep the text quickly accessible.\\

\textbf{Notation:} For a small category $\cC$, we shall write $\cC_{0}$
and $\cC_{1}$ for its sets of objects and morphisms. If $F\from
\cC\to\cD$ is a functor, then $F_{0}$ and $F_{1}$ denotes the maps on
objects and morphisms. If $F,G \from \cC\to\cD$ are functors and
$\alpha_{x}\from F(x)\to G(x)$ is natural, then we also write $\alpha
\from F\Rightarrow G$ if we want to emphasise thinking of $\alpha$ as a
map $\alpha\from \cC_{0}\to \cD_{1}$. Moreover, we write $\Delta_{\cC}$
(or shortly $\Delta$ if $\cC $ is understood) for the diagonal embedding
$\cC\to\cC\times\cC$. Unless stated otherwise, all categories are
assumed to be small.

\section{Principal 2-bundles}

In this section we introduce principal 2-bundles. For those readers who
wonder what the 2 refers to: throughout this paper we are working
\textit{in} the 2-category of categories, where objects are given by
categories, morphisms by functors between categories and 2-morphisms by
natural transformations between functors. It is \textit{not} the case that
the things we call 2-something are 2-categories by themselves (just as a
set is not a category but an object in the category of sets).

Although the repeated term ``\dots such that there exist natural
equivalences\dots '' is quite annoying we state it explicitly every
time it occurs, because it is the source of the additional structure
in the theory (compared to ordinary bundle theory), which deserves to
be pointed out. On the other hand, the occurring coherence conditions
in terms of commutative diagrams
can safely be neglected at first reading, because they will all be
trivially satisfied later on.

\begin{definition}\label{def:2-group}
 A \textit{weak 2-group} is a monoidal category, in which each morphism
 is invertible and each object is weakly invertible. We spell this out
 for convenience. It is given by a category $\cG$ together with a
 multiplication functor $\otimes \from\cG\times\cG\to\cG$ (mostly
 written as $g\cdot h:=\otimes (g,h)$), an object $\idObject$ of
 $\cG$ and natural equivalences
 \begin{equation*}
  \alpha_{g,h,k}\from (g\cdot h)\cdot k\to g\cdot(h\cdot k),
 \end{equation*}
 $l_{g}\from \idObject\cdot g\to g$ and
 $r_{g}\from g\cdot\idObject\to g$, such that the diagrams {\small
 \[
  \xymatrix@=10pt{ &&{\scriptstyle (g\cdot h)\cdot (k\cdot l)}
  \ar[drr]^(.5){\displaystyle\alpha_{g,h,k\cdot l}}\\
  {\scriptstyle ((g\cdot h)\cdot k)\cdot
  l}\ar[urr]^(.5){\displaystyle\alpha_{g\cdot
  h,k,l}}\ar[ddr]_{\displaystyle\alpha_{g,h,k}\cdot \id_{l}}&&&&
  {\scriptstyle g\cdot(h\cdot (k\cdot l))}\\
  \\
  &{\scriptstyle (g\cdot (h\cdot k))\cdot
  l}\ar[rr]_{\displaystyle\alpha_{g,h\cdot k,l}}&&{\scriptstyle
  g\cdot((h\cdot k)\cdot l)}\ar[uur]_{\displaystyle\id_{g}\cdot
  \alpha(h,k,l)} }
 \]
 } and {\small
 \begin{equation*}
  \xymatrix{
  {\scriptstyle (g\cdot \idObject)\cdot h}\ar[dr]_(.4){\displaystyle r_{g}\cdot\id_{h}}\ar[rr]^{\displaystyle \alpha_{g,\idObject,h}}&& {\scriptstyle g\cdot (\idObject \cdot h)}\ar[dl]^(.4){\displaystyle \id_{g}\cdot l_{h}}\\
  &  {\scriptstyle g\cdot h}}
 \end{equation*}
 }commute. Moreover, we require that each morphism is invertible and
 that for each object $x$ there exists an object $\ol{x}$ such that
 $x\cdot \ol{x}$ and $\ol{x}\cdot x$ are isomorphic to $\idObject$.

 A \textit{2-group} is a weak 2-group, together with a \emph{coherent
 choice} of (weak) inverses, given by an additional functor
 $\iota \from \cG\to\cG$ and natural equivalences
 $i_{g}\from g\cdot \iota(g)\to \idObject$ and
 $e_{g}\from \iota(g)\cdot g\to \idObject$, such that {\small
 \begin{equation*}
  \xymatrix{
  {\scriptstyle(g\cdot \iota(g))\cdot g}\ar[d]_{\displaystyle \alpha_{g,\iota(g),g}}\ar[rr]^(.55){\displaystyle i_{g}\cdot\id_{g}}&&{\scriptstyle\idObject\cdot g}\ar[r]^(.6){\displaystyle l_{g}}&\ar[d]^{\displaystyle \id_{g}}{\scriptstyle g}\\
  {\scriptstyle g\cdot (\iota(g)\cdot g)}\ar[rr]^(.55){\displaystyle \id_{g}\cdot e_{g}}&&{\scriptstyle g\cdot \idObject}\ar[r]^(.6){\displaystyle r_{g}}&{\scriptstyle g}
  }
 \end{equation*}
 }commutes. Morphisms of 2-groups are defined to be weakly monoidal
 functors of the underlying monoidal category (cf.\
 \cite{baezLauda04Higher-DimensionalAlgebraV:2-Groups}).
\end{definition}

Our main reference for 2-groups is
\cite{baezLauda04Higher-DimensionalAlgebraV:2-Groups}, where our 2-groups
are called coherent 2-groups. As also mentioned in
\cite{baezLauda04Higher-DimensionalAlgebraV:2-Groups}, this is what is
also called a (coherent) category with group structure (cf.\
\cite{laplaza83CoherenceForCategoriesWithGroupStructures:AnAlternativeApproach},
\cite{ulbrich81KohaerenzInKategorienMitGruppenstruktur}).

\begin{example}
(cf.\ \cite[Ex.\ 34]{baezLauda04Higher-DimensionalAlgebraV:2-Groups})
Let $\cC$ be a category and $\cAut_{w}(\cC)$ be the category of
equivalences of $\cC$. Then $\cAut_{w}(\cC)$ is a weak 2-group with
respect to composition of functors and natural transformations. There
is also a coherent version of this 2-group, cf.\ \cite[Ex.\
35]{baezLauda04Higher-DimensionalAlgebraV:2-Groups}.
\end{example}

\begin{example}
 Let $G$ be a group, $A$ be abelian and $f\from G\times G\times G\to A$
 be a group  cocycle, i.e., we have
 \begin{equation}\label{eqn:groupCocycle}
  f(gh,k,l)+f(g,h,kl)=f(g,h,k)+f(g,hk,l)+f(h,k,l)
 \end{equation}
 for $g,h,k,l\in G$. Then we define a category $\cG_{f}$ by setting
 $\Obj(\cG_{f}):=G$ and
 \begin{equation*}
  \Hom(g,g')=\left\{\begin{array}{ll}A &\tx{ if }g=g' \\ \emptyset & \tx{ else}\end{array}\right.
 \end{equation*}
 with the composition coming from group multiplication in $A$. Then
 \begin{equation*}
  (g\xrightarrow{a}g)\cdot (h\xrightarrow{b} h):=gh\xrightarrow{ab}gh
 \end{equation*}
 defines a multiplication functor on $\cG_{f}$ and
 \begin{equation*}
  (g,h,k)\mapsto\left( g\cdot h\cdot k\xrightarrow{f(g,h,k)}g\cdot h\cdot k\right)
 \end{equation*}
 defines a a natural equivalence. That this natural equivalences make
 the diagrams from Definition \ref{def:2-group} commutes is equivalent
 to \eqref{eqn:groupCocycle}. Thus, $\cG_{f}$ is a weak 2-group and each
 2-group is equivalent to such a 2-group (cf.\ \cite[Sect.\
 8.3]{baezLauda04Higher-DimensionalAlgebraV:2-Groups}).
\end{example}

\begin{definition}
 A \emph{smooth 2-space} is a small category $\cM$ such that $\cM_{0}$,
 $\cM_{1}$ and ${\cM_{1}}_{\,s\!\!}\times_{t} \cM_{1}$ are smooth
 manifolds and all structure maps are smooth.

 A \emph{smooth functor} $F\from \cM\to \cM'$ between smooth 2-spaces is
 a functor such that $F_{0}$ and $F_{1}$ are smooth maps. Likewise, a
 \emph{smooth natural transformation} $\alpha\from F\Rightarrow G$
 between smooth functors is a natural transformation which is smooth as
 a map $\cM_{0}\to \cM_{1}$.

 Eventually, a \emph{smooth equivalence} between smooth 2-spaces $\cM$
 and $\cM'$ is a smooth functor $F\from \cM\to\cM'$ such that there
 exist a smooth functor $G\from\cM'\to\cM$ and smooth natural
 equivalences $G\circ F\Rightarrow \id_{\cM}$ and
 $F\circ G\Rightarrow \id_{\cM'}$.
\end{definition}

In our context, a manifold refers to a Hausdorff space, which is locally
homeomorphic to open subsets of a locally convex space such that the
coordinate changes are smooth (cf.\ Appendix
\ref{app:differentialCalculusOnSpacesOfMappings}). This definition of a
smooth 2-space is a bit more rigid than the concept used frequently in
the literature for it requires each space really to be a smooth manifold
and not just a smooth (or diffeological) space (cf.\
\cite{baezSchreiber05HigherGaugeTheory},
\cite{BaezHoffnung08Convenient-Categories-of-Smooth-Spaces}). For our
present aim, this concept of a smooth 2-space suffices.
Most of the smooth 2-spaces that appear in this article are in fact Lie groupoids,
but there is no need to restrict to Lie groupoids a priori.

\begin{example}
The easiest example of a smooth 2-space is simply a smooth manifold as
space of objects with only identity morphisms and the obvious
structure maps.
\end{example}

\begin{definition}
 A (strong) \textit{Lie 2-group} is a 2-group which is a smooth 2-space
 at the same time, such that the functors and natural equivalences
 occurring in the definition of a 2-group are smooth.
\end{definition}

In general, it is quite restrictive to require all functors and natural
equivalences to be smooth (cf.\ \cite[Sect.\
2]{Wockel08Categorified-central-extensions-etale-Lie-2-groups-and%
-Lies-Third-Theorem-for-locally-exponential-Lie-algebras} and
\cite{Henriques08Integrating-Lsb-infty-algebras}). However, the major
part of this paper deals with strict 2-groups, where the definition is
appropriate (cf.\ \cite[Sect.\
2]{Wockel08Categorified-central-extensions-etale-Lie-2-groups-and%
-Lies-Third-Theorem-for-locally-exponential-Lie-algebras}). 

We now consider how smooth 2-groups may act on smooth 2-spaces.

\begin{definition}\label{def:G-2-space}
 Let $\cG$ be a Lie 2-group. Then a \textit{smooth $\cG$-2-space} is a
 smooth 2-space $\cM$ together with a smooth action, i.e., a smooth
 functor $\rho\from \cM\times\cG \to \cM$ (mostly written as
 $x.g:=\rho(x,g)$) and smooth natural equivalences
 $\nu  \from \rho \circ (\rho \times \id_{\cG})\Rightarrow \rho \circ (\id_{\cM}\times \otimes ) $
 and
 $\xi  \from \rho \circ (\id_{\cM}\times \idObject)\Rightarrow \id_{\cM}$
 such that the diagrams {\small
 \[
  \xymatrix@=10pt{ &&{\scriptstyle (x.g).(h\cdot k)}
  \ar[drr]^(.5){\displaystyle \nu _{x,g,h\cdot k}}\\
  {\scriptstyle ((x.g).h).k}\ar[urr]^(.5){\displaystyle\nu_{x.g,h,k}}
  \ar[ddr]_{\displaystyle\nu_{x,g,h}.\id_{k}}&&&&{\scriptstyle x.(g\cdot
  (h\cdot k))} \\
  \\
  &{\scriptstyle (x.(g\cdot h))\cdot
  k}\ar[rr]_{\displaystyle\nu_{x,g\cdot h,k}}&&{\scriptstyle x.((g\cdot
  h)\cdot k)}\ar[uur]_{\displaystyle\id_{x}. \alpha(g,h,k)} }
 \]
 } and {\small
 \begin{equation*}
  \xymatrix{
  {\scriptstyle (x. \idObject). g}\ar[dr]_(.4){\displaystyle \xi_{x}.\id_{g}}\ar[rr]^{\displaystyle \nu_{x,\idObject,g}}&& {\scriptstyle x.(\idObject \cdot g)}\ar[dl]^(.4){\displaystyle \id_{x}. l_{g}}\\
  &  {\scriptstyle x.g}}
 \end{equation*}
 }commute. A \textit{morphism} between $\cG$-2-spaces is a smooth
 functor $F\from \cM\to\cM'$, and a smooth natural equivalence
 $\sigma^{F} \from F\circ \rho\Rightarrow \rho '\circ (F\times \id_{\cG})$
 such that {\small
 \[
  \xymatrix@=10pt{ &&{\scriptstyle F(x).(g\cdot h)}
  \ar[drr]^(.5){\displaystyle \sigma^{F}_{x,g\cdot h}}\\
  {\scriptstyle (F(x).g).h }\ar[urr]^(.5){\displaystyle \nu'_{x,g,h}
  }\ar[ddr]_{\displaystyle \sigma^{F}_{x,g}.\id_{h}}&&&& {\scriptstyle
  F(x.(g\cdot h))}\\
  \\
  &{\scriptstyle F(x.g).h }\ar[rr]_{\displaystyle
  \sigma^{F}_{x.g,h}}&&{\scriptstyle F((x.g).h)}\ar[uur]_{\displaystyle
  F(\nu_{x,g,h})} }
 \]
 }commutes. A \textit{2-morphism} between two morphisms
 $F,F'\from \cM\to\cM'$ of $\cG$-2-spaces is a smooth natural
 transformation $\tau \from F\Rightarrow F'$, such that the diagram
 {\small
 \begin{equation*}
  \xymatrix{
  F(x).g \ar[d]_{\sigma^{F}_{x,g}} \ar[r]^(.45){\tau_{x}.\id_{g}}&
  F'(x).g\ar[d]^{\sigma^{F'}_{x,g}}\\
  F(x.g)\ar[r]_(.5){\tau_{x.g}}& F'(x.g)
  }
 \end{equation*}
 }commutes. An \emph{equivalence} of smooth $\cG$-2-spaces is a morphism
 $F\from \cM\to\cM'$ such that there exists a morphism
 $F'\from\cM'\to\cM$ and 2-morphisms $F \circ F'\Rightarrow \id_{\cM'}$
 and $F' \circ F\Rightarrow \id_{\cM}$. In this case, $F'$ is called a
 \emph{weak inverse} of $F$.
\end{definition}

Since we shall only consider smooth actions of Lie 2-groups on smooth
2-spaces we suppress this adjective in the sequel. We are now
ready to define principal bundles, whose base is a 2-space with only
identity morphisms.

\begin{definition}\label{def:principal-2-bundle}
 Let $\cG$ be a Lie 2-group and $M$ be a smooth manifold (viewed as a
 smooth 2-space with only identity morphisms). A \emph{principal
 $\cG$-2-bundle} over $M$ is a locally trivial $\cG$-2-space over $M$.
 More precisely, it is a smooth $\cG$-2-space $\cP$, together with a
 smooth functor $\pi\from\cP\to M$, such that there exist an open cover
 $(U_{i})_{i\in I}$ of $M$ and equivalences
 $\Phi_{i}\from\left.\cP\right|_{U_{i}}\to U_{i}\times\cG$ of
 $\cG$-2-spaces (where $\cG$ acts on $U_{i}\times\cG$ by right
 multiplication on the second factor). Moreover, we require
 $\left.\pi \right|_{U_{i}}=\pr_{1}\circ \Phi_{i}$ and
 $\pi \circ \ol{\Phi}_{i}=\pr_{1}$  on the nose for a
 weak inverse $\ol{\Phi}_{i}$ of $\Phi_{i}$. Various diagrams, emerging
 from the natural equivalences, are required to commute to ensure
 coherence (cf.\ \cite[Sect.\
 2.5]{bartels06HigherGaugeTheoryI:2-Bundles}).

 A \emph{morphism} of principal $\cG$-2-bundles over $M$ is a morphism
 $\Phi\from \cP\to \cP'$ of $\cG$-2-spaces satisfying
 $\pi' \circ \Phi=\pi$, and a 2-morphism between two morphisms of
 principal $\cG$-2-bundles is a 2-morphism of the underlying
 morphisms of strict $\cG$-2-spaces.
 As above, various diagrams are required to commute (cf.\ \cite[Sect.\
 2.5]{bartels06HigherGaugeTheoryI:2-Bundles}).
\end{definition}

We suppress an explicit statement of the coherence conditions for
brevity. We do not need them in the sequel, as we shall restrict to
cases where most natural equivalences are required to be the identity
transformation. Note that in the previous sense, a principal $\cG$-2-bundle is
``locally trivial'', i.e., each $\left.\cP\right|_{U_{i}}$ is
equivalent to $U_{i}\times \cG$. In particular, each principal $\cG$-2-bundle
is a Lie groupoid.

\begin{lemma}
Principal $\cG$-2-bundles over $M$, together with their morphisms
and smooth natural equivalences between morphisms form a 2-category
\mbox{2-$\cBun(M,\cG)$}.
\end{lemma}

\begin{proof}
 It is easily checked, that \mbox{2-$\cBun(M,\cG)$} actually is a
 sub-2-category of the 2-category of small categories, functors and natural
 transformations.
\end{proof}

\section{Classification of principal 2-bundles by \v{C}ech cohomology}

So far, we have clarified the categorification procedure for principal
bundles. We now stick to more specific examples of
these bundles, which are classified by non-abelian \v{C}ech
cohomology. The idea is to strictify everything that concerns the
action in case of a strict structure group. Strictification means for us
to require natural transformations to be the identity.

\begin{definition}
 A \emph{strict 2-group} is a 2-group $\cG$, where all natural
 equivalences between functors occurring in the definition of a 2-group
 are the identity. A morphism of strict 2-groups is a weak monoidal
 functor $F\from \cG\to\cG'$ with $F(g\cdot h)=F(G)\cdot F(h)$.
\end{definition}

We promised to keep the reader away from 2-categories. However, many
formulae and calculations become intuitively understandable in a
graphical representation, which we shortly outline in the following
remark. The reader who wants to neglect this representation may do so,
we shall provide at each stage explicit formulae.

\begin{remark}\label{rem:2-groupAs2-category}
 For a diagrammatic interpretation of various formulae and arguments, it
 is convenient to view a strict 2-group $\cG$ not only as a category,
 but also as a 2-category with one object. The reader unfamiliar with
 2-categories may understand this as the association of an arrow
 \begin{equation*}
  \xymatrix{\bullet \ar[r]^{g}&\bullet}
 \end{equation*}
 between one fixed object $\bullet$, which we assign to each object $g$
 of $\cG$, and the association of a 2-arrow
 \begin{equation*}
  \xymatrix{
  \bullet
  \ar@/^1pc/[rr]^{g}_{\ }="s"
  \ar@/_1pc/[rr]_{g'}^{\ }="t"
  &&
  \bullet
  \ar@{=>}^{h} "s"; "t"
  }
 \end{equation*}
 between the arrows $g$ and $g'$, which we assign to each morphism
 $h\from g\to g'$ in $\cG$. Then composition in $\cG$ is depicted by the
 vertical composition
 \begin{equation*}
  \xymatrix{
  \bullet
  \ar@/^1.6pc/[rr]^{g}_{\ }="a"
  \ar[rr]|{\stackrel{{\scriptstyle g'}}{\ }}="b"
  \ar@/_1.6pc/[rr]_{g''}^{\ }="c"
  &&
  \bullet
  \ar@{=>}^{h} "a"; "b"
  \ar@{=>}^{h'} "b"; "c"
  }=
  \xymatrix{
  \bullet
  \ar@/^1.6pc/[rr]^{g}_{\ }="a"
  \ar@/_1.6pc/[rr]_{g''}^{\ }="b"
  &&
  \bullet
  \ar@{=>}^{h' \circ h} "a"; "b"
  }
 \end{equation*}
 of 2-arrows. These diagrams should cause no confusion with the kind of
 diagrams from the previous section, where objects were represented by
 points and morphisms were represented by arrows. The latter kind of
 diagrams will not occur any more in the sequel.

 The multiplication functor on objects is then depicted by the
 horizontal concatenation
 \begin{equation*}
  \xymatrix{\bullet \ar[r]^{g}&\bullet \ar[r]^{\ol{g}}&\bullet}=
  \xymatrix{\bullet \ar[rr]^{g\cdot \ol{g}}&&\bullet}
 \end{equation*}
 of arrows. On morphisms, multiplication is depicted by the horizontal
 concatenation
 \begin{equation*}
  \xymatrix{
  \bullet
  \ar@/^1pc/[rr]^{g}_{\ }="s1"
  \ar@/_1pc/[rr]_{g'}^{\ }="t1"
  &&
  \bullet
  \ar@/^1pc/[rr]^{\ol{g}}_{\ }="s2"
  \ar@/_1pc/[rr]_{\ol{g}'}^{\ }="t2"
  &&
  \bullet
  \ar@{=>}^{h} "s1"; "t1"
  \ar@{=>}^{\ol{h}} "s2"; "t2"
  }= \xymatrix{
  \bullet
  \ar@/^1pc/[rrr]^{g\cdot \ol{g}}_{\ }="s"
  \ar@/_1pc/[rrr]_{g'\cdot \ol{g}'}^{\ }="t"
  &&&
  \bullet
  \ar@{=>}^{h\cdot \ol{h}} "s"; "t"
  }
 \end{equation*}
 of 2-arrows.
\end{remark}

We shall only deal with strict 2-groups in the following text, and
there are many different ways to describe them. In \cite[p.\
3]{baezLauda04Higher-DimensionalAlgebraV:2-Groups} one finds the
following list (which is explained in detail in
\cite{Porst08Strict-2-Groups-are-Crossed-Modules}). A strict
2-group is
\begin{itemize}\label{list:differentDescriptionOf2-groups}
\item a strict monoidal category in which all objects and morphisms
are invertible,
\item a strict 2-category with one object in which all 1-morphisms and
2-morphisms are invertible,
\item a group object in categories (also called a \textit{categorical
group}),
\item a category object in groups (also called \textit{internal
category} in groups),
\item a crossed module.
\end{itemize}

For the present text, the interpretation of a strict 2-group as a
crossed module (and vice versa) will be of central interest, so we
recall this concept and relate it to 2-groups.

\begin{definition}
A \textit{crossed module} consists of two groups $G,H$, an action
$\alpha \from G\to \Aut(H)$ and a homomorphism $\beta \from H\to G$
satisfying
\begin{align}
\beta (\alpha (g).h)&=g\cdot \beta (h)\cdot 
g^{-1}\label{eqn:crossedModule_equivariance}\\
\alpha (\beta (h)).h'&=h\cdot h'\cdot
h^{-1}\label{eqn:crossedModule_peiffer}.
\end{align}
\end{definition}

\begin{remark}\label{rem:2-groupsFromCrossedModules}
 From a crossed module, one can build a 2-group as follows, cf.\
 \cite{Porst08Strict-2-Groups-are-Crossed-Modules}, \cite[Prop.\
 16]{bartels06HigherGaugeTheoryI:2-Bundles},
 \cite{forrest-barker02GroupObjectsAndInternalCategories} and
 \cite[Sect.\
 XII.8]{Mac-Lane98Categories-for-the-working-mathematician}. The set of
 objects is $G$ and the set of morphisms is $H\rtimes G$. Each element
 $(h,g)\in H\rtimes G$ defines a morphism from $g$ to $\beta(h)g$. This
 is graphically represented by a 2-arrow
 
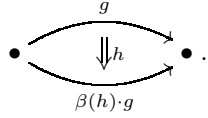
\begin{equation*}
  \xymatrix{
  \bullet
  \ar@/^1pc/[rr]^{g}_{\ }="s"
  \ar@/_1pc/[rr]_{\beta(h)\cdot g}^{\ }="t"
  &&
  \bullet\; .
  \ar@{=>}^{h} "s"; "t"
  }
 \end{equation*}
 Composition in the category is given by
 $(h',\beta (h)\cdot g)\circ (h,g):=(h'\cdot h,g)$. This is depicted by
 defining
 \begin{equation}\label{eqn:2-groupsFromCrossedModules1}
  \xymatrix{
  \bullet
  \ar@/^1.6pc/[rr]^{g}_{\ }="a"
  \ar[rr]|{\stackrel{{\scriptstyle \beta(h)\cdot g}}{\ }}="b"
  \ar@/_1.6pc/[rr]_{\beta(h')\cdot \beta (h)\cdot g}^{\ }="c"
  &&
  \bullet \phantom{.}
  \ar@{=>}^{h} "a"; "b"
  \ar@{=>}^{h'} "b"; "c"
  }:=
  \xymatrix{
  \bullet
  \ar@/^1.6pc/[rr]^{g}_{\ }="a"
  \ar@/_1.6pc/[rr]_{\beta(h'\cdot h)\cdot g}^{\ }="b"
  &&
  \bullet\; .
  \ar@{=>}^{h'\cdot h} "a"; "b"
  }
 \end{equation}
 Consequently, $(e,g)$ defines the identity of $g$. One easily checks
 that the space of composable pairs of morphisms is
 $H\rtimes(H\rtimes G)$, where $H\rtimes G$ acts on $H$ by
 $(h,g).h'=\alpha (\beta (h)\cdot g).h'$ and on this space, composition
 is given by the homomorphism $(h',(h,g))\mapsto (h'\cdot h,g)$.
 Similarly one shows that the space of composable triples of morphisms
 is $H\rtimes (H\rtimes (H\rtimes G))$ and the associativity in $H$
 yields the associativity of composition. The multiplication functor is
 determined by the group multiplication in $H\rtimes G$. On objects, it
 is given by multiplication in $G$ and is  thus depicted by
 \begin{equation*}
  \xymatrix{\bullet \ar[r]^{g}&\bullet \ar[r]^{\ol{g}}&\bullet}:=
  \xymatrix{\bullet \ar[rr]^{g\cdot \ol{g}}&&\bullet\quad .}
 \end{equation*}
 On morphisms, it is given by multiplication in $H\rtimes G$ and thus
 depicted by
 \begin{equation*}
  \xymatrix{
  \bullet
  \ar@/^1pc/[rr]^{g}_{\ }="s1"
  \ar@/_1pc/[rr]_{\beta(h)\cdot g}^{\ }="t1"
  &&
  \bullet
  \ar@/^1pc/[rr]^{\ol{g}}_{\ }="s2"
  \ar@/_1pc/[rr]_{\beta(\ol{h})\cdot \ol{g}}^{\ }="t2"
  &&
  \bullet
  \ar@{=>}^{h} "s1"; "t1"
  \ar@{=>}^{\ol{h}} "s2"; "t2"
  }:= \xymatrix{
  \bullet
  \ar@/^1pc/[rrr]^{g\cdot \ol{g}}_{\ }="s"
  \ar@/_1pc/[rrr]_{\beta(h)\cdot  g\cdot\beta(\ol{h}) \cdot \ol{g}}^{\ }="t"
  &&&
  \bullet
  \ar@{=>}^{h\cdot (g.\ol{h})} "s"; "t"
  }
 \end{equation*}
 (note that the target of $(h\cdot (g.\ol{h}),g\cdot\ol{g})$ is in fact
 $\beta(h)\cdot  g\cdot\beta(\ol{h}) \cdot \ol{g}$ by
 \eqref{eqn:crossedModule_equivariance}). Likewise, the inversion
 functor is determined by inversion in $H\rtimes G$. All this together
 defines a strict 2-group.

 The reverse construction is also straightforward. For a strict 2-group
 one checks that objects and morphisms are groups themselves and that
 all structure maps are group homomorphisms. Then one sets $H$ to be the
 kernel of the source map and $G$ to be the space of objects. Then $G$
 acts on $H$ by $g.h=\id_{g}\cdot h\cdot \id_{g^{-1}}$ and $\beta$ is
 given by the restriction of the target map to $H$.

 Historically, crossed modules arose first in the work of Whitehead on
 homotopy 2-types
 \cite{Whitehead46Note-on-a-previous-paper-entitled-On-adding-relations-to-homotopy-groups},
 and the equivalence of crossed modules and 2-groups was established by
 Brown and Spencer in
 \cite{BrownSpencer76G-groupoids-crossed-modules-and-the-fundamental-groupoid-of-a-topological-group}
 and by Loday in
 \cite{Loday82Spaces-with-finitely-many-nontrivial-homotopy-groups}. A
 detailed exposition of the equivalence of the 2-categories of strict
 2-groups and of crossed modules is given in
 \cite{Porst08Strict-2-Groups-are-Crossed-Modules}.
\end{remark}

\begin{remark}\label{rem:whiskering}
 To put more involved diagrams as the ones occurring in the previous
 remark into formulae, one has to replace
 \begin{equation}\label{eqn:whiskeringFromTheLeft}
  \vcenter{\xymatrix{
  \bullet \ar[r]^{\ol{g}}&\bullet
  \ar@/^1pc/[r]^{g}_{\ }="s"
  \ar@/_1pc/[r]_{\beta(h)\cdot g}^{\ }="t"
  &
  \bullet
  \ar@{=>}^{h} "s"; "t"
  }}\quad\quad\tx{ by }\quad\quad
  \vcenter{\xymatrix{
  \bullet
  \ar@/^1pc/[r]^{\ol{g}\cdot g}_(.3){\ }="s"
  \ar@/_1pc/[r]_{\ol{g}\cdot \beta(h)\cdot g}^(.3){\ }="t"
  &
  \bullet
  \ar@{=>}^{\ol{g}.h} "s"; "t"
  }}
 \end{equation}
 on a $\bullet$ with at least one incoming two outgoing 1-arrows and on
 a $\bullet$ with at least one outgoing and two incoming 1-arrows one
 has to replace
 \begin{equation}\label{eqn:whiskeringFromTheRight}
  \vcenter{\xymatrix{
  \bullet
  \ar@/^1pc/[r]^{g}_{\ }="s"
  \ar@/_1pc/[r]_{\beta(h)\cdot g}^{\ }="t"
  &
  \bullet \ar[r]^{\ol{g}}&\bullet
  \ar@{=>}^{h} "s"; "t"
  }}\quad\quad\tx{ by }\quad\quad
  \vcenter{\xymatrix{
  \bullet
  \ar@/^1pc/[r]^{g\cdot \ol{g}}_{\ }="s"
  \ar@/_1pc/[r]_{\beta(h)\cdot g\cdot \ol{g}}^{\ }="t"
  &
  \bullet .
  \ar@{=>}^{h} "s"; "t"
  }}
 \end{equation}
 This procedure is called \emph{whiskering} and corresponds exactly to
 multiplication of $(h,g)\in H\rtimes G$ with $(e,\bar{g})$ from the
 left (in \eqref{eqn:whiskeringFromTheLeft}) and from the right (in
 \eqref{eqn:whiskeringFromTheRight}). By performing these substitutions
 one ends up in a diagram which depicts compositions of 2-arrows, which
 can be performed as in \eqref{eqn:2-groupsFromCrossedModules1}. That
 various ways of doing these substitutions yield the same result is
 encoded in the axioms of a crossed module. We shall see how this works
 in practise in Definition \ref{def:non-abelianCechCohomology}.
\end{remark}

\begin{definition}
A \textit{smooth crossed module} is a
crossed module $(G,H,\alpha ,\beta)$ such that $G,H$ are
Lie groups, $\beta$ is smooth and the map $G\times H\to H$,
$(g,h)\mapsto \alpha (g).h$ is smooth.
\end{definition}

\begin{remark}
Given a smooth crossed module, the pull-back of composable pairs of
morphism is $H\rtimes(H\rtimes G)$ and the space of composable triples
of morphisms is $H\rtimes(H\rtimes (H\rtimes G))$. From the
description in Remark \ref{rem:2-groupsFromCrossedModules} it follows
that all structure maps are smooth and thus the corresponding 2-group
actually is a Lie 2-group.
\end{remark}

For a strict Lie 2-group to define a smooth crossed module it is
necessary that the kernel of the source map is a Lie subgroup (cf.\
Example \ref{rem:2-groupsFromCrossedModules}). In finite dimensions
this is always true but an infinite-dimensional Lie group may possess
closed subgroups that are no Lie groups (cf.\ \cite[Ex.\
III.8.2]{bourbakiLie}). So in the differentiable setup we take smooth
crossed modules as the concept from the list on page
\pageref{list:differentDescriptionOf2-groups} representing all its equivalent
descriptions.

\begin{definition}
 If $\cG$ is a strict Lie 2-group, then a $\cG$-2-space is called
 \textit{strict} if all natural equivalences between functors in the
 definition of a $\cG$-2-space are identities. Likewise, a
 \textit{morphism} between strict $\cG$-2-spaces is a morphism
 $F\from\cM\to\cM'$ of the underlying $\cG$-2-spaces with $\sigma^{F}$
 the identity transformation.
 A \emph{2-morphism} between two morphisms $F,F'\from \cM\to\cM'$ of
 strict $\cG$-2-spaces is 2-morphism
 $\tau\from F\Rightarrow F'$, satisfying $\tau(x.g)=\tau(x).\id_{g}$. An
 \emph{equivalence} of strict $\cG$-2-spaces is a morphism
 $F\from \cM\to\cM'$ such that there exist a morphism
 $F'\from\cM'\to\cM$ and 2-morphisms (all of strict $\cG$-2-spaces)
 $F \circ F'\Rightarrow \id_{\cM'}$ and
 $F' \circ F\Rightarrow \id_{\cM}$. In this case, $F'$ is called a
 \emph{weak inverse} of $F$.
\end{definition}

The crucial point for the following text shall be that we restrict to
strict $\cG$-2-spaces, i.e., we do \emph{not} allow non-identical
natural equivalences to occur in any axiom that concerns the action
functor $\rho \from\cM\times\cG\to\cM$.

\begin{definition}
 If $\cG$ is a strict Lie 2-group, then we call a principal $\cG$-2-bundle
 $\cP$ \textit{semi-strict} if $\cG$ acts strictly on $\cP$ and the
 local trivialisations may be chosen to be equivalences of strict
 $\cG$-2-spaces. A morphism between semi-strict principal
 $\cG$-2-bundles $\cP$ and $\cP'$ over $M$ (or a \textit{semi-strict
 bundle morphism} for short) is a morphism $\Phi\from \cP\to\cP'$ of
 strict $\cG$-2-spaces satisfying $\pi'\circ \Phi=\pi$. Likewise, a
 2-morphism between two morphisms of semi-strict principal
 $\cG$-2-bundles is a 2-morphism between the underlying morphisms of
 strict $\cG$-2-spaces.
\end{definition}

Note that our concept of semi-strictness differs from the one used in
\cite{bartels06HigherGaugeTheoryI:2-Bundles}, which is a normalisation
requirement on the cocycles classifying principal 2-bundles. A
\textit{strict} principal $\cG$-2-bundle would require that all natural
equivalences occurring in its definition can be chosen to be the
identity. However, this definition is slightly too rigid for a treatment
of non-abelian cohomology. On the other hand, many generalisations are
possible by increasingly admitting various additional natural
equivalences to be non-trivial.

\begin{lemma}
Semi-strict principal $\cG$-2-bundles over $M$, together with their
morphisms and 2-morphisms form a
2-category \mbox{$ss$-2-$\cBun(M,\cG)$}.
\end{lemma}

It shall turn out that semi-strict principal 2-bundles are the immediate
generalisations of principal bundles, as long as one is interested in
geometric structures defined over ordinary (smooth) spaces, classified in
terms of \v{C}ech cohomology. 
Moreover, semi-strict principal 2-bundles 
are also linked to gerbes as follows.

\begin{remark}\label{rem:from-principal-2-bundles-to-gerbes}
 In general, a gerbe over $X$ is a locally transitive and locally non-empty stack
 in groupoids (cf.\
 \cite{moerdijk02IntroductionToTheLanguageOfStacksAndGerbes}). This
 means that if $U\mapsto F(U)$ is the underlying fibered category of the
 stack, then one assumes $X=\bigcup \{U|F(U)\neq \emptyset\}$ and that for objects
 $a,b$ of $F(U)$,
 each $x\in U$ has an open neighbourhood $V$ in $U$ with
 at least one arrow $\left.a\right|_{V}\to \left.b\right|_{V}$.

 Now if $\cP$ is a principal $\cG$-2-bundle with structure group coming
 from the smooth crossed module $H\to\Aut(H)$ (and assume that $\Aut(H)$
 is a Lie group, for instance if $\pi_{0}(H)$ is finitely generated,
 cf.\ \cite{bourbakiLie}), then this defines a gerbe by
 \begin{equation*}
  U\mapsto  \{\tx{sections of }\left.\cP\right|_{U}\}
 \end{equation*}
 for $U$ an open subset of $M$ (cf.\ \cite[Thm.\
 3.1]{moerdijk02IntroductionToTheLanguageOfStacksAndGerbes}).
\end{remark}

From now on we shall assume that $\cG$ is a smooth Lie 2-group
arising from a smooth crossed module $(\alpha ,\beta ,H,G)$. We spend
the rest of this section on the classification problem for semi-strict
principal $\cG$-2-bundles over $M$.

It shall turn out that principal $\cG$-2-bundles over $M$ are classified
(in an appropriate sense) by the non-abelian cohomology $\check
H(M,\cG)$. There are many treatments of non-abelian cohomology in the
literature, i.e., \cite{bartels06HigherGaugeTheoryI:2-Bundles},
\cite{breenMessing05DifferentialGeometryOfGerbes},
\cite{breen92OnTheClassificationOf2-GerbesAnd2-Stacks},
\cite{giraud71CohomologieNonAbelienne},
\cite{dedecker60SurLaCohomologieNonAbelienne}, and our definition is
essentially the same (with the usual minor conventional differences).
Note that we did not put a degree (or dimension) to $\check H(M,\cG)$,
for the kind of 2-group that one takes for $\cG$ determines its degree
in ordinary \v{C}ech-cohomology. For instance, $\check
H(M,\cG)=\check{H}^{2}(M,A)$ if $\cG$ is associated to the crossed
module $A\to\{*\}$ (for $A$ abelian) and $\check
H(M,\cG)=\check{H}^{1}(M,G)$ if $\cG$ is associated to the crossed
module $\{*\}\to G$ (for $G$ arbitrary).

\begin{definition}\label{def:non-abelianCechCohomology}
 An element in $\check{H}(M,\cG)$ is represented by an open cover
 $(U_{i})_{i\in I}$ of $M$, together with a collection of smooth maps
 $g_{ij}\from U_{ij}\to G$ and $h_{ijk}\from U_{ijk}\to H$ (where
 multiple lower indices refer to multiple intersections). These maps are
 required to satisfy point-wise
 $\beta(h_{ijk})\cdot g_{ij}\cdot g_{jk}=  g_{ik}$ on $U_{ijk}$, i.e.,
 \begin{equation}\label{eqn:cocyc1}
  \vcenter{\xymatrix{ 
  &\bullet\ar[dr]^{g_{jk}}\ar@{=>}[d]^(.66){h_{ijk}}&\\ 
  \bullet \ar[rr]_{g_{ik}}\ar[ur]^{g_{ij}}&&\bullet
  }}
 \end{equation}
 and $h_{ikl} \cdot h_{ijk}=h_{ijl}\cdot (g_{ij}.h_{jkl})$ on $U_{ijkl}$,
 i.e.,
 \begin{equation}\label{eqn:cocyc2}
  \vcenter{\xymatrix@=4em{
  \bullet \ar[r]^{g_{jk}}&\bullet \ar[d]^{g_{kl}}\\
  \bullet \ar[u]^{g_{ij}} \ar[ur]|{\blabel{g_{ik}}}="a"
  \ar[r]_{g_{il}}^(.4){\ }="b"&\bullet
  \ar@{=>}^(.6){h_{ijk}} "1,1"+/dr 1em/; "a"
  \ar@{=>}^{h_{ikl}} "1,2"+/d:(1,-0.6) 2em/; "b"
  }}
  =
  \vcenter{\xymatrix@=4em{
  \bullet \ar[r]^{g_{jk}}\ar[dr]|{\blabel{g_{jl}}}="a"&
  \bullet \ar[d]^{g_{kl}}\\
  \bullet \ar[u]^{g_{ij}}  \ar[r]_{g_{il}}^(.6){\ }="b"&\bullet
  \ar@{=>}_(.6){h_{jkl}} "1,2"; "a"
  \ar@{=>}_{h_{ijl}} "1,1"+/d:(1,0.6) 2em/; "b"
  }}
 \end{equation}
 (note that the occurrence of $g_{ij}$ in the formula is caused by
 whiskering, cf.\ Remark \ref{rem:whiskering}). We furthermore require
 $g_{ii}=e_{G}$ and $h_{ijj}=h_{jji}=e_{H}$ point-wise. We call such a
 collection $(U_{i},g_{ij},h_{ijk})$ a (non-abelian, normalised,
 $\cG$-valued) \emph{cocycle} on $M$.

 Two cocycles $(U_{i},g_{ij},h_{ijk})$ and
 $(U_{i'}',g'_{i'j'},h'_{i'j'k'})$ are called \emph{cohomologous} (or \emph{equivalent}) if there
 exist a common refinement $(V_{i})_{i\in J}$ (i.e.,
 $V_{i}\se U_{\bar i}$ and $V_{i}\se U'_{\bar i'}$ for each $i\in J$ and
 some $\bar i\in I$ and $\bar i'\in I'$) and a collection of smooth maps
 $\gamma_{i}\from V_{i}\to G$ and $\eta_{ij}\from V_{ij}\to H$,
 satisfying point-wise
 $\gamma_{i}\cdot g'_{\bar i'\bar j'}=\beta(\eta_{ij})\cdot g_{\bar i\bar j}\cdot \gamma_{j}$
 on $V_{ij}$, i.e.,
 \begin{equation}\label{eqn:cobound1}
  \vcenter{\xymatrix{
  \bullet \ar[r]^{g_{\bar i\bar j}}_{\ }="s"
  \ar[d]_{\gamma_{i}} & \bullet \ar[d]^{\gamma_{j}}\\
  \bullet \ar[r]_{g'_{\bar i'\bar j'}}^{\ }="t" & \bullet 
  \ar@{=>}^{\eta_{ij}} "1,2"; "2,1"
  }}
 \end{equation}
 and
 $\eta_{ik}\cdot h_{\bar i\bar j\bar k}=\gamma_{i}. h'_{\bar i'\bar j'\bar k'}\cdot \eta_{ij}\cdot g_{\bar i\bar j}.\eta_{jk}$
 on $V_{ijk}$, i.e.,
 \begin{equation}\label{eqn:cobound2}
  \vcenter{\xymatrix@=2.25em{
  &\bullet \ar[dr]^{g_{\bar j\bar k}}\\
  \bullet \ar[dd]_{\gamma_{i}}
  \ar[rr]|{\blabel{g_{\bar i\bar k}}}^{\ }="t1"_{\ }="s2" 
  \ar[ur]^{g_{\bar i\bar j}}&&
  \bullet\ar[dd]^{\gamma_{k}}\\
  \\
  \bullet \ar[rr]_{g'_{\bar i'\bar k'}}^{\ }="t2"&& \bullet
  \ar@{=>}^(.7){h_{\bar i\bar j\bar k}} "1,2"; "t1"
  \ar@{=>}^{\eta_{ik}} "2,3"+/dl 1.25em/; "4,1"+/ur 1.25em/
  }}
  =
  \vcenter{\xymatrix@=2.25em{
  &\bullet \ar[dr]^{g_{\bar j\bar k}}|{\vphantom{\blabel{g}}}="s2" 
  \ar[dd]|{\blabel{\gamma_{j}}}="t4"\\
  \bullet \ar[dd]_{\gamma_{i}} 
  \ar[ur]^{g_{\bar i\bar j}}|{\vphantom{\blabel{g}}}="s1"&&
  \bullet\ar[dd]^{\gamma_{k}}\\
  & \bullet \ar[dr]|{\blabel{g'_{\bar j'\bar k'}}}="t2"\\
  \bullet \ar[ur]|{\blabel{g'_{\bar i'\bar j'}}}="t1"
  \ar[rr]_{g'_{\bar i'\bar k'}}^{\ }="t3"&& 
  \bullet\; .
  \ar@{=>}^{\eta_{ij}} "1,2"+/d:(3,-1) 1.5em/; "4,1"+/d:(-3,1) 1.5em/
  \ar@{=>}^(.6){\eta_{jk}} "2,3"; "3,2"
  \ar@{=>}^(.8){h'_{\bar i'\bar j'\bar k'}} "3,2"; "t3"
  }}
 \end{equation}
 We furthermore require $\eta_{ii}=e_{H}$ point-wise. Sometimes,
 $(V_{i},\gamma_{i},\eta_{ij})$ is also called a \emph{coboundary}
 between $(U_{i},g_{ij},h_{ijk})$ and $(U'_{i'},g'_{i'j'},h'_{i'j'k'})$.
 It is easily checked by taking refinements and point-wise products in
 $\cG$ that this defines in fact an equivalence relation and we denote
 by $\check{H}(M,\cG)$ the resulting set of equivalence classes of
 cocycles.
\end{definition}

Note that our normalisation conditions $g_{ii}=e_{G}$,
$h_{ijj}=h_{iij}=e_{H}$ and $\eta_{ii}=e_{H}$ do \emph{not} imply
$g_{ij}=g_{ji}^{-1}$, $h_{ijk}=h_{jik}^{-1}$ and $\eta_{ij}=h_{ji}^{-1}$
as one might expect. Note also that one obtains the non-abelian cohomology as used in the
texts, mentioned above, when one takes for a (connected) Lie group $G$ the crossed module
$G\to\Aut(G)$, induced by the conjugation homomorphism.

\begin{remark}
 The previous definition is not arbitrary but is the natural
 generalisation of the following idea. If $G$ is a Lie group, then an
 ordinary $G$-valued cocycle on $M$ is given by an open cover
 $(U_{i})_{i\in I}$ and smooth $g_{ij}\from U_{ij}\to G$ satisfying
 $g_{ij}g_{jk}=g_{jk}$ and $g_{ii}=e_{G}$ point-wise. But this is the
 same as a smooth functor from the \v{C}ech groupoid
 of $(U_{i})_{i\in I}$ to the smooth one-object
 category $BG$. Likewise, coboundaries between cocycles
 are given by natural transformations between the
 corresponding functors on refined covers.

 If $\cG$ is a Lie 2-group, then we can view $\cG$ as a 2-category
 $B\cG$ with only one object (cf.\ Remark
 \ref{rem:2-groupAs2-category}). Moreover
 the \v{C}ech groupoid can also be viewed as a
 2-category with only identity 2-morphisms. Then the cocyles arise as
 pseudo (or weak) 2-functors from
 this 2-category to $B\cG$
 and coboundaries as pseudonatural transformations between them on
 refined covers
 (cf.\ \cite[Sect.\ 7.5]{borceux94HandbookOfCategoricalAlgebraVolume1} or \cite[App.\ A]{SchreiberWaldorf08Connections-on-non-abelian-Gerbes-and-their-H%
 olonomy} for the terminology).
\end{remark}

The point of view from the previous remark also shows that cocycles
actually have more categorical structure, i.e., there is a 2-category of
cocycles, coming from the 2-category structure of pseudofunctors,
pseudonatural transformations and modifications from \v{C}ech groupoids
to $\cG$. However, this lies beyond the scope of the present article and
we shall not elaborate on this additional structure.\\

We shall now show how to obtain a cocycle from a principal 2-bundle.

\begin{remark}\label{rem:fromBundlesToCocycles}
 The argument explained below reoccurs frequently in the following
 construction. For $U\se M$, each morphism of strict $\cG$-2-spaces (or
 strictly equivariant functor) $\Psi\from U\times \cG\to U\times \cG$,
 which is the identity in the first component, is given by
 $(x,g)\mapsto (x,g(x)^{-1}\cdot g)$ on objects for a map
 $g\from U\to G$ and
 $(x,(h,g))\mapsto (x,(g(x)^{-1}.h,g(x)^{-1}\cdot g))$ on morphisms (it
 is determined by its values on the subcategory
 $(U\times \{e_{G}\},U\times (e_{H},e_{G}))$ and the artificial
 inversions is taken in order to match our other conventions, in
 particular \eqref{eqn:fromBundlesToCocycles_1}). If $\Psi$ is smooth,
 then $g$ is smooth and vice versa. For two different such
 $\Psi_{1},\Psi_{2}$, a 2-morphism $\Psi_{1}\Rightarrow \Psi_{2}$
 between morphisms of strict $\cG$-2-spaces is given by
 $(x,g)\mapsto (x,(g_{1}(x)^{-1}.h(x)^{-1},g_{1}(x)^{-1}\cdot g))$ for a
 unique map $h\from U\to H$ satisfying
 \begin{equation}\label{eqn:fromBundlesToCocycles_1}
  g_{2}=(\beta\circ h)\cdot g_{1}.
 \end{equation}
 Moreover, each map $h\from U\to H$, satisfying
 \eqref{eqn:fromBundlesToCocycles_1}, defines a natural equivalence
 $\Psi_{1}\Rightarrow \Psi_{2}$ which is smooth if and only if $h$ is
 smooth.

 For a principal $\cG$-2-bundle $\cP$ we now construct a cocycle
 $z(\cP)$ as follows. We choose an open cover $(U_{i})_{i\in I}$ and
 local trivialisations
 $\Phi_{i}\from \left.\cP\right|_{U_{i}}\to U_{i}\times\cG$, as well as
 weak inverses
 $\ol{\Phi}_{i}\from U_{i}\times\cG\to \left.\cP\right|_{U_{i}}$ and
 2-morphisms
 $\tau_{i}\from \Phi_{i}\circ\ol{\Phi}_{i}\Rightarrow \id_{U_{i}\times \cG}$
 and
 $\ol{\tau}_{i}\from \ol{\Phi}_{i}\circ{\Phi}_{i}\Rightarrow \id_{\left.\cP\right|_{U_{i}}}$.
 For each pair $i,j\in I$, we consider the composition
 $\Phi_{j}\circ \ol{\Phi}_{i}\from U_{ij}\times \cG\to U_{ij}\times \cG$,
 of local trivialisations. This is of the above form and thus determined
 by a smooth map $g_{ij}\from U_{ij}\to G$, i.e.,
 \begin{equation*}
  \vcenter{\xymatrix{
  & \left.\cP\right|_{U_{ij}} \ar[dr]^{\Phi_{j}} & \\
  U_{ij}\times \cG \ar[ur]^{\ol{\Phi}_{i}} \ar[rr]_{g_{ij}}^{\ }="t"&& 
  U_{ij}\times \cG\;,
  \ar@{=} "1,2"+/d 2em/; "t"+/u 1em/
  }}
 \end{equation*}
 or, equivalently, $\Phi_{j}(\ol{\Phi}_{i}(x,e))=(x,g_{ij}^{-1}(x))$.
 Substituting $\Phi_{i}$ by $(g_{ii}^{-1} \circ \pi)\cdot \Phi_{i}$, we
 may assume that $\Phi_{i} \circ\ol{\Phi}_{i}=\id_{U_{i}\times\cG}$ on
 the nose and thus $g_{ii}=e_{G}$.

 For each $i_{1},\dots ,i_{n}\in I$ with $n\geq 2$ we define the functor
 \[
  \Psi_{i_{1}\dots i_{n}}:=\Phi_{i_{n}}\circ \ol{\Phi}_{i_{n-1}}\circ
  \Phi_{i_{n-1}}\circ\dots \circ \ol{\Phi}_{i_{2}} \circ
  \Phi_{i_{2}}\circ \ol{\Phi }_{i_{1}}.
 \]
 On objects, $\Psi_{ijk}$ is then given by
 $(x,g)\mapsto (x,((g_{ij}(x)\cdot g_{jk}(x))^{-1}\cdot g))$. With
 $\ol{\tau}_{j}\from \ol{\Phi}_{j}\circ \Phi_{j}\Rightarrow \id_{U_{j}\times\cG}$,
 composition of natural transformations yields a smooth equivalence
 $\Psi_{ijk}\Rightarrow \Psi_{ik}$ given by a smooth map
 $ h_{ijk}\from U_{ijk}\to H$ satisfying
 $(\beta \circ h_{ijk})\cdot g_{ij}\cdot g_{jk}= g_{ik}$ by
 \eqref{eqn:fromBundlesToCocycles_1}, i.e.,
 \begin{equation*}
  \vcenter{\xymatrix{
  & U_{ijk}\times\cG \ar[dr]^{g_{jk}}\\
  U_{ijk}\times\cG\ar[rr]_{g_{ik}}^{\ }="t" \ar[ur]^{g_{ij}}&&U_{ijk}\times\cG
  \ar@{=>}^(.6){h_{ijk}} "1,2";"t"
  }}:= \hspace{-2em}
  \vcenter{\xymatrix@=2em{
  & U_{ijk}\times\cG \ar[dddr]^{g_{jk}}_{\ }="jk"
  \ar@/^1pc/[dd]|{\ol{\Phi}_{j}}\\ \\
  & \left.\cP\right|_{U_{ijk}} \ar[dr]|{\Phi_{k}} 
  \ar@/^1pc/[uu]|{{\Phi}_{j}}\\
  U_{ijk}\times \cG \ar[ur]|{\ol{\Phi}_{i}} 
  \ar[uuur]^{g_{ij}}_{\ }="ij"
  \ar[rr]_{g_{ik}}^{\ }="ik" && U_{ijk}\times\cG \; .
  \ar@{=>}|{\blabel{\ol{\tau}_{j}}}"1,2";"3,2"
  \ar@{=}"3,2";"ij"
  \ar@{=}"3,2";"jk"
  \ar@{=}"3,2"+/d 1.8em/; "ik"+/u .8em/
  }} 
 \end{equation*}
 Due to \eqref{eqn:crossedModule_peiffer}, we thus have
 \begin{equation}\label{eqn:fromBundlesToCocycles_2}
 \Phi_{k}(\ol{\tau}_{j}(\ol{\Phi}_{i}(x,e)))= (g_{ik}(x)^{-1}. h_{ijk}(x)^{-1},(g_{ij}(x)\cdot g_{jk}(x))^{-1}).
 \end{equation} 
 From $\Phi_{i} \circ \ol{\Phi}_{i}=\id_{U_{i}\times \cG}$ it follows
 that the above natural transformation is the identity if two
 neighbouring indices agree and thus that $h_{iij}=h_{ijj}=e_{H}$.
 Moreover, composition of natural transformations yields smooth
 equivalences $\Psi_{ijkl}\Rightarrow \Psi_{ijl}$ and
 $\Psi_{ijkl}\Rightarrow \Psi_{ikl}$, and the compositions
 $\Psi_{ijkl}\Rightarrow \Psi_{ijl}\Rightarrow \Psi_{il}$ and
 $\Psi_{ijkl}\Rightarrow \Psi_{ikl}\Rightarrow \Psi_{il}$ coincide (cf.\
 \cite[Prop.\ 1.3.5]{borceux94HandbookOfCategoricalAlgebraVolume1}).
 Spelling this out leads to
 $h_{ikl}\cdot h_{ijk}=h_{ijl}\cdot (g_{ij}.h_{jkl})$. The corresponding
 diagram can be obtained from plugging the diagram, defining $h_{ijk}$,
 into \eqref{eqn:cocyc2}.

 The cocycle $z(\cP)$ depends on the choices of $(U_{i})_{i\in I}$,
 $\Phi_{i}$, $\ol{\Phi}_{i}$ and $\ol{\tau}_{i}$. However, a special
 case of the following argument shows that two different such choices
 lead to cohomologous cocycles.

 Assume for the moment that we have fixed some choice of the previous
 data for each $\cP$. Given a morphism $\Phi\from \cP\to\cP'$ of
 semi-strict principal $\cG$-2-bundles over $M$, we shall show that the
 cocycles $z(\cP)$ and $z(\cP')$, constructed from that data, are
 cohomologous. First, we choose a common refinement $(V_{i})_{i\in J}$ of
 $(U_{i})_{i\in I}$ and $(U'_{i'})_{i'\in I'}$, i.e.,
 $V_{i}\se U_{\bar i}$ and $V_{i}\se U'_{\bar i'}$ for each $i\in J$ and
 some $\bar i\in I$ and $\bar i'\in I'$. Clearly, all the choices of
 $\Phi_{\bar i}$'s and $\tau_{\bar i}$'s restrict to $V_{i}$.

 For each $i\in J$, we consider the morphism
 $\Phi_{\bar i'}'\circ \Phi \circ \ol{\Phi}_{\bar i}\from V_{i}\times \cG\to V_{i}\times \cG$.
 Since $\pi'\circ \Phi=\pi$, this is of the form described initially and
 thus determined by a smooth map $\gamma_{i}\from V_{i}\to G$, i.e.,
 \begin{equation*}
  \vcenter{\xymatrix{
  & \left.\cP\right|_{V_{i}}\ar[r]^{\Phi}_{\ }="s"
  &\left.\cP'\right|_{V_{i}}
  \ar[dr]^{\Phi'_{\bar i'}}\\
  V_{i}\times \cG \ar[ur]^{\ol{\Phi}_{\bar i}}\ar[rrr]_{\gamma_{i}}^{\ }="t"
  &&& V_{i}\times\cG.
  \ar@{=}"s"+/d 1.25em/; "t"+/u .75em/
  }}
 \end{equation*}
 As above, $\ol{\tau}_{j}$ and ${\ol{\tau}'_{\bar i}}^{-1}$ give rise to
 a smooth natural equivalence
 \begin{equation*}
  \underbrace{\Phi'_{\bar j'}\circ \Phi \circ \ol{\Phi}_{\bar j}}_{=\gamma_{j}}
  \circ \underbrace{\Phi_{\bar j}\circ \ol{\Phi}_{\bar i}}_{=g_{\bar i\bar j}}
  \Rightarrow
  \underbrace{\Phi'_{\bar j'}\circ \ol{\Phi}'_{\bar i'}}_{=g'_{\bar i'\bar j'}}
  \circ 
  \underbrace{{\Phi}'_{\bar i'}\circ
  \Phi\circ\ol{\Phi}_{\bar i}}_{=\gamma_{i}}
 \end{equation*}
 given by a smooth map $\eta_{ij}\from V_{ij}\to H$ satisfying
 $\beta(\eta_{ij})\cdot g_{\bar i\bar j}\cdot \gamma_{j}=\gamma_{i}\cdot g'_{\bar i'\bar j'}$
 by \eqref{eqn:fromBundlesToCocycles_1}, i.e.,
 \begin{equation*}
  \vcenter{\xymatrix{
  V_{ij}\times \cG \ar[r]^{\gamma_{j}}_{\ }="s"
  & V_{ij}\times \cG  \\
  V_{ij}\times \cG  \ar[u]^{g_{\bar i\bar j}}\ar[r]_{\gamma_{i}}^{\ }="t" &
  V_{ij}\times \cG 
  \ar[u]_{g'_{\bar i'\bar j'}}
  \ar@{=>}^{\eta_{ij}} "1,1"; "2,2"
  }}:=
  \vcenter{\xymatrix{
  V_{ij}\times \cG \ar[rrr]^{\gamma _{j}} _(.49){\ }="a"
  \ar@/^1pc/[dr]|{\blabel{\ol{\Phi}_{\bar j}}} 
  &&&
  V_{ij}\times\cG\\
  & \left.\cP\right|_{V_{ij}}  \ar[r]|{\blabel{\Phi}} ="cent2"
  \ar@/^1pc/[ul]|{\blabel{\Phi_{\bar j}}}
  & \left.\cP'\right|_{V_{ij}} \ar[ur]|{\Phi_{\bar j'}'} 
  \ar@/_1pc/[dr]|{\blabel{\Phi_{\bar i'}'}}\\
  V_{ij}\times\cG \ar[rrr]_{\gamma_{i}}^(.49){\ }="c" \ar[uu]^{g_{\bar i\bar j}}_{\ }="d"
  \ar[ur]|{\blabel{\ol{\Phi}_{\bar i}}}&&& V_{ij} \times\cG\; . \ar@/_1pc/[ul]|{\ol{\Phi}_{\bar i'}'} \ar[uu]_{g'_{\bar i'\bar j'}}^{\ }="b"
  \ar@{=} "a" +/d 1.25em/; "cent2" +/u 1.25em/
  \ar@{=} "cent2" +/d 1.25em/; "c" +/u 1.25em/
  \ar@{=} "2,2" +/l 2.75em/; "d" +/r 1em/
  \ar@{=} "b" +/l 1em/; "2,3" +/r 3.25em/
  \ar@{=>} "1,1"; "2,2"|{\blabel{\ol{\tau}_{\bar j}}}
  \ar@{=>} "2,3"; "3,4"|{\blabel{{\ol{\tau}_{\bar i'}'}^{-1}}}
  }}
 \end{equation*}
 Now the compositions
 \begin{equation*}
  \Phi_{\bar k'}'\circ \Phi \circ \ol{\Phi}_{\bar k}\circ \Phi_{\bar k}\circ
  \ol{\Phi}_{\bar j}\circ
  \Phi_{\bar j}\circ \ol{\Phi}_{\bar i} 
  \xymatrix{\ar@{=>}[r]^{h_{\bar i\bar j\bar k}}&}
  \Phi_{\bar k'}'\circ \Phi \circ \ol{\Phi}_{\bar k}\circ \Phi_{\bar k}\circ
  \ol{\Phi}_{\bar i} 
  \xymatrix{\ar@{=>}[r]^{\eta_{ik}}&}
  \Phi_{\bar k'}'\circ \ol{\Phi}_{\bar i'}' \circ \Phi_{\bar i'}' 
  \circ \Phi \circ
  \ol{\Phi}_{\bar i} 
 \end{equation*}
 
 and
 \begin{multline*}
  \Phi_{\bar k'}'\circ \Phi \circ \ol{\Phi}_{\bar k}\circ \Phi_{\bar k}\circ
  \ol{\Phi}_{\bar j}\circ
  \Phi_{\bar j}\circ \ol{\Phi}_{\bar i} 
  \xymatrix{\ar@{=>}[r]^{\eta_{jk}}&}
  \Phi_{\bar k'}'\circ \ol{\Phi}_{\bar j'}'\circ \Phi_{\bar j'}'\circ \Phi \circ
  \ol{\Phi}_{\bar j}\circ
  \Phi_{\bar j}\circ \ol{\Phi}_{\bar i} \\
  \xymatrix{\ar@{=>}[r]^{\eta_{ij}}&}
  \Phi_{\bar k'}'\circ \ol{\Phi}_{\bar j'}'\circ \Phi_{\bar j'}'\circ
  \ol{\Phi}_{\bar i'}'\circ
  \Phi_{\bar i'}'\circ  \Phi \circ \ol{\Phi}_{\bar i} 
  \xymatrix{\ar@{=>}[r]^{h'_{\bar i'\bar j'\bar k'}}&}
  \Phi_{\bar k'}'\circ \ol{\Phi}_{\bar i'}' \circ \Phi_{\bar i'}' 
  \circ \Phi \circ
  \ol{\Phi}_{\bar i}  
 \end{multline*}
 coincide (cf.\ \cite[Prop.\
 1.3.5]{borceux94HandbookOfCategoricalAlgebraVolume1}). Spelling this
 out leads to
 \[
  \eta_{ik}\cdot h_{\bar i\bar j\bar k}=\gamma_{i}. h'_{\bar i'\bar
  j'\bar k'}\cdot \eta_{ij}\cdot g_{\bar i\bar j}.\eta_{jk}
 \]
 The corresponding diagram can be obtained from plugging the diagrams,
 defining $h_{ijk}$ and $\eta _{ij}$ into \eqref{eqn:cobound2}.

 From the previous argument it also follows that different choices of
 $(U_{i})_{i\in I}$, $\Phi_{i}$, $\ol{\Phi}_{i}$, $\tau_{i}$ and
 $\ol{\tau}_{i}$ in the construction of $z(\cP)$ lead to cohomologous
 cocycles. In fact, if we apply the construction to the identity of
 $\cP$, then the resulting $\gamma_{i}\from U_{i}\to G$ and
 $\eta_{ij}\from U_{ij}\to H$ yield the desired coboundary.
\end{remark}

\begin{proposition}\label{prop:fromBundlesToCocycles}
The construction from the previous remark assigns to each semi-strict principal $\cG$-2-bundle $\cP$ an element $[z(\cP)]\in\check{H}(M,\cG)$. 
Moreover, if for $\cP$ and $\cP'$ there exists a semi-strict bundle morphism
$\Phi\from\cP\to\cP'$, then $[z(\cP)]=[z(\cP')]$.
\end{proposition}

What is special about the construction in Remark
\ref{rem:fromBundlesToCocycles} is that we can construct a principal
2-bundle $\cP_{z}$ out of a given 2-cocycle $z=(U_{i},g_{ij},h_{ijk})$
with $[z]=[z(\cP)]$. To illustrate the construction we first recall that
for a given group valued cocycle $g_{ij}\from U_{ij}\to G$ a
corresponding principal bundle is given by \[ P=\bigcup_{i\in
I}\,\,\{i\}\times U_{i}\times G/\sim, \] where $(i,x,g_{ij}(x)\cdot
g)\sim (j,x, g)$. The main idea is to modify this construction by
introducing a refined identification. Thinking categorically, this means
that we do not only remember that objects are isomorphic (equivalent),
but also track consistently the different isomorphisms that may exist.
Identifying isomorphic objects in $\cP_{z}$ then leads to the
construction below of the underlying principal $G/\beta (H)$-bundle (cf.\
Corollary \ref{rem:BandOfAGerbe}). Moreover, the following construction
shall be tailored to generalise the fact that for an ordinary principal
bundle, the fibres of the projection map are equivalent to the structure
group $G$ as right $G$-spaces.

\begin{remark}\label{rem:fromCocyclesToBundles}
 Given a 2-cocycle $z=(U_{i},g_{ij},h_{ijk})$, we define the category
 $\cP_{z}$ by
 \begin{align*}
  \Obj(\cP_{c})&=\bigcup_{i\in I}\,\,\{i\}\times U_{i}\times G\\
  \Mor(\cP_{c})&=\bigcup_{i,j\in I}\{(i,j)\}\times U_{ij}\times H\times
  G
 \end{align*}
 with the obvious smooth structure. It shall be clear in the sequel that
 all maps are smooth with respect to this structure, so we shall not
 comment on this any more. We want the structure maps to make the
 identification of
 \begin{equation*}
  (i,x,g)\xrightarrow{(i,j,x,h,g)}(j,x,g')\quad\tx{ with }\quad
  \vcenter{\xymatrix{
  \bullet \ar[rr]^{g}_{\ }="s" \ar[dr]_{g_{ij}}&&\bullet \\
  & \bullet \ar[ur]_{g'}
  \ar@{=>}^{h} "s"; "2,2" +/u 1em/
  }}
 \end{equation*}
 (the latter diagram is in $\cG$) an equivalence of categories. We thus
 require $g'=g_{ij}^{-1}\cdot\beta(h)\cdot g$. Consequently, the space
 of composable morphisms is
 $\coprod U_{ijk}\times H\times G$, which we
 also endow with the obvious smooth structure. Then composition in
 $\cP_{z}$ is defined by setting the composition
 \begin{equation*}
  (i,x,g)
  \xrightarrow{(i,j,x,h,g)}
  (j,x,\underbrace{g_{ij}^{-1}\cdot\beta(h)\cdot g}_{=:g'})
  \xrightarrow{(j,k,x,h',g')}
  (k,x,g_{jk}^{-1}\cdot\beta(h')\cdot g')
 \end{equation*}
 to be
 \begin{equation*}
  (i,x,g)\xrightarrow{(i,k,x,h_{ijk}\cdot (g_{ij}.h') \cdot
  h,g)}(k,x,g_{ik}^{-1}\cdot \beta(h_{ijk}\cdot g_{ij}.h' \cdot h)\cdot g),
 \end{equation*}
 i.e., it is induced by the composition
 \begin{equation}\label{eqn:fromCocyclesToBundles1}
  \vcenter
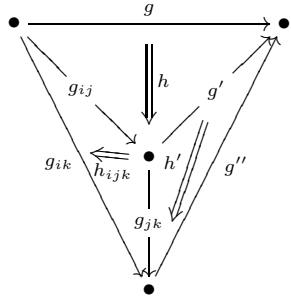
{\xymatrix@=4em{
  \bullet \ar[dr]|{\blabel{g_{ij}}} \ar[rr]^{g}_{\ }="s1" 
  \ar[ddr]_{g_{ik}}^{\ }="t2" && \bullet\\
  & \bullet \ar[ru]|{\blabel{g'}}="t1" \ar[d]|{\blabel{g_{jk}}} \\
  & \bullet \ar[uur]_{g''}
  \ar@{=>}^(.45){h} "s1"+/d .5em/; "2,2"+/u 1.25em/
  \ar@{=>}_(.45){h'} "t1"+/d:(3,-1) 1.25em/; "3,2"+/d:(-3,1) 2.75em/
  \ar@{=>}^(.35){{h_{ijk}}} "2,2"+/l .75em/; "t2"+/l .125em/
  }}
 \end{equation}
 in $\cG$. That this definition satisfies source-target matching follows
 from $\beta(h_{ijk})\cdot g_{ij}\cdot g_{jk}=g_{ik}$. Moreover,
 $(i,i,x,e,g)$ defines the identity of $(i,x,g)$ and
 $(i,j,x,h,g)^{-1}=(j,i,x,g_{ij}^{-1}.((h\cdot h_{iji})^{-1}),g_{ij}^{-1}\cdot\beta(h)\cdot g)$.
 That composition is associative follows from
 $h_{ikl} \cdot h_{ijk}=h_{ijl}\cdot g_{ij}.h_{jkl}$ (a corresponding
 equality of diagrams may be obtained from sticking together
 \eqref{eqn:fromCocyclesToBundles1} and \eqref{eqn:cocyc2}). Thus we
 obtain a smooth 2-space $\cP_{z}$ with an obvious morphism
 $\cP_{z}\to M$.

 The right action of $\cG$ on $\cP_{c}$ is given by
 $(i,x,g).\ol{g}:= (i,x,g\cdot \ol{g})$ on objects and by
 \[
  (i,j,x,h,g).(\ol{h},\ol{g}):=(i,j,x,h \cdot (%
g.\ol{h}),
  g\cdot \ol{g})
 \]
 on morphisms, i.e., it is induced by the horizontal composition
 \begin{equation*}
  \vcenter{\xymatrix{
  \bullet \ar@/^1pc/[rr]^{{g}}_{\ }="s1" 
  \ar@/_1pc/[rr]_{g_{ij}\cdot g'}^{\ }="t1" &&
  \bullet \ar@/^1pc/[rr]^{\ol{g}}_{\ }="s2" 
  \ar@/_1pc/[rr]_{\ol{g}'}^{\ }="t2" && \bullet
  \ar@{=>}^{{h}} "s1"; "t1"
  \ar@{=>}^{\ol{h}} "s2"; "t2"
  }}
  =
  \vcenter{\xymatrix{
  \bullet \ar@/^1pc/[rrr]^{g\cdot \ol{g}}_(.4){\ }="s" 
  \ar@/_1pc/[rrr]_{g_{ij}\cdot g'\cdot\ol{g}'}^(.4){\ }="t" &&& \bullet \quad .
  \ar@{=>}^{h \cdot(g. \ol{h})} "s"; "t"
  }}
 \end{equation*}
 
 We want the local trivialisations $\ol{\Phi}_{i}$ to be given by the
 canonical inclusion $U_{i}\times \cG\to \cP_{z}$ and $\Phi_{i}$ to be
 tailored such that the natural equivalences
 $\ol{\tau}_{i}\from \ol{\Phi}_{i}\circ \Phi_{i}\Rightarrow \id_{\left.\cP_{z}\right|_{U_{i}}}$
 are given by
 \begin{equation*}
  (j,x,g)\mapsto \left(\xymatrix{(i,x,g_{ij}\cdot g) 
  \ar[rr]^{(i,j,x,e,g_{ij}\cdot g)}&& (j,x,g)}\right).
 \end{equation*}
 We thus set
 $ \Phi_{i}\from \left.\cP_{z}\right|_{U_{i}}\to U_{i}\times \cG $ to be
 induced by the assignment
 \begin{equation*}
  \vcenter{\xymatrix{
  (i,x,g_{ij} g) 
  \ar[d]_{(i,j,x,e,g_{ij} g)} 
  \ar[rr]^(.4){(i,i,x,h^{*},g_{ij} g)} &&
  (i,x,g_{ik}g_{jk}^{-1}\beta(h) g)\\
  (j,x,g) 
  \ar[rr]^(.45){(j,k,x,h,g)} && 
  (k,x,g_{jk}^{-1} \beta (h) g)
  \ar[u]^{(i,k,x,e,g_{ik} g_{jk}^{-1} \beta(h)g)^{-1}}
  }} \mapsto 
  \vcenter{\xymatrix{
  \bullet \ar@/^1pc/[rrr]^{(x,g_{ji}^{-1}g)}_(.3){\ }="s"
  \ar@/_1pc/[rrr]_{(x,(g_{jk}g_{ki})^{-1}g)}^(.3){\ }="t" &&& \bullet
  \ar@{=>}^{(x,h^{*},g_{ji}^{-1}g)} "s"; "t"
  }}
 \end{equation*}
 with $h^{*}=h_{ijk}\cdot g_{ij}.h$, i.e., $\Phi_{i}$ is defined by
 \begin{align*}
  (j,x,g) & \mapsto (x,g_{ij}\cdot g)  \\
  (j,k,x,h,g) & \mapsto (x,h_{ijk}\cdot g_{ij}.h,g_{ij}\cdot g).
 \end{align*}
 In fact, $\Phi_{i}\circ \ol{\Phi}_{i}=\id_{U_{i}\times \cG}$ on the
 nose and
 $ \ol{\tau}_{i}\from\id_{\left.\cP_{c}\right|_{U_{i}}}\Rightarrow \ol{\Phi}_{i}\circ\Phi_{i}$
 is then given by $(j,x,g)\mapsto (i,j,x,e,g_{ij}g)$.
\end{remark}

\begin{proposition}\label{prop:fromCocyclesToBundles}
 For each $\cG$-valued cocycle $z$ on $M$, the principal $\cG$-2-bundle
 $\cP_{z}$ over $M$, constructed in the previous remark, has
 $[z(\cP_{z})]=[z]$. Moreover, if $z$ is equivalent to $z'$, then there
 exists a semi-strict principal $\cG$-2-bundle $\cP$ over $M$, and two
 semi-strict bundle morphisms $\Phi\from\cP\to \cP_{z}$ and
 $\Phi'\from\cP\to \cP_{z'}$.
\end{proposition}

\begin{proof}
 Applying the construction of $z(\cP_{z})$ from Remark
 \ref{rem:fromBundlesToCocycles} to the bundle $\cP_{z}$, constructed in
 Remark \ref{rem:fromCocyclesToBundles} (and choosing in this
 construction $(U_{i})_{i\in I}$, $\Phi_{i}$, $\ol{\Phi}_{i}$,
 $\tau_{i}$ and $\ol{\tau}_{i}$ as defined in Remark
 \ref{rem:fromCocyclesToBundles}) shows that with these choices we have
 in fact $z(\cP_{z})=z$.

 In order to verify the second claim we first consider the case of two
 cohomologous cocycles $z=(U_{i},g_{ij},h_{ijk})$ and
 $z'=(U'_{i'},g'_{i'j'},h_{i'j'k'})$, where $(U_{i})_{i\in I}$ is a
 refinement of $(U'_{i'})_{i'\in I'}$. I.e., we have
 $U_{i}\se U'_{\bar i'}$ for some $\bar i\in I'$, and the coboundary is
 given by smooth maps $\gamma_{i}\from U_{i}\to G$ and
 $\eta_{ij}\from U_{ij}\to H$. In this case we set $\cP$ to be
 $\cP_{z}$, $\Phi$ to be $\id_{\cP_{z}}$, and it remains to construct
 $\Phi'\from\cP_{z}\to\cP_{z'}$ as follows. First, \eqref{eqn:cobound1}
 induces
 \begin{equation*}
  \vcenter{\xymatrix@=2em{
  & \bullet \ar[dr]^{g_{ij}}\\
  \bullet \ar[ur]^{\gamma_{i}^{-1}}\ar[dr]_{g'_{\bar i\bar j}} & & \bullet\\
  & \bullet \ar[ur]_{\gamma_{j}^{-1}}
  \ar@{=>}|{\blabel{\gamma_{i}^{-1}.\eta_{ij}}} "1,2"+/d 1em/; "3,2"+/u 1em/
  }}:=
  \vcenter{\xymatrix@=2em{
  && \bullet \ar[dr]^{\gamma_{j}}\\
  \bullet \ar[r]^{\gamma_{i}^{-1}} & \bullet 
  \ar[ur]^{g_{ij}}\ar[dr]_{\gamma_{i}} & & 
  \bullet \ar[r]^{\gamma_{j}^{-1}} & \bullet\\
  && \bullet \ar[ur]_{g'_{\bar i\bar j}}
  \ar@{=>}|{\blabel{\eta_{ij}}} "1,3"+/d 1em/; "3,3"+/u 1em/
  }}
 \end{equation*}
 and this induces $\Phi'$ by the assignment
 \begin{equation}\label{eqn:fromCocyclesToBundles2}
  \vcenter{\xymatrix@C=4em{
  \bullet \ar[rr]^{g\vphantom{\gamma_{i}^{-1}}}_{\ }="s" 
  \ar[dr]_{g_{ij}}&&\bullet \\
  & \bullet \ar[ur]_{g'} \ar@{}[r]_{\vphantom{\gamma_{j}^{-1}}}&
  \ar@{=>}^{h} "s"; "2,2" +/u 1em/
  }}
  \mapsto
  \vcenter{\xymatrix@C=4em{
  \bullet \ar[r]^{\gamma_{i}^{-1}} \ar[dr]_{g_{\bar i\bar j}'}& 
  \bullet \ar[rr]^{g}_{\ }="s" \ar[dr]|{\blabel{g_{ij}}}
  &
  &\bullet \\
  & \bullet \ar[r]_{\gamma_{j}^{-1}}
  & \bullet \ar[ur]_{g'} \ar@{}[r]_{\vphantom{\gamma_{j}^{-1}}}&
  \ar@{=>}^{h} "s"; "2,3" +/u 1em/
  \ar@{=>}|{\gamma_{i}^{-1}.\eta_{ij}} "1,2"; "2,2"
  }}
 \end{equation}
 that is
 \begin{align*}
  \Phi'(i,x,g)&=(\bar i,x,\gamma_{i}^{-1}\cdot g)\\
  \Phi'(i,j,x,h,g)&=(\bar i,\bar j,x,\eta_{ij}\cdot(\gamma_{i}^{-1}. h),\gamma_{i}^{-1}\cdot g).
 \end{align*}
 That $\Phi'$ satisfies source-target matching follows from
 $\gamma_{i}\cdot g'_{\bar i'\bar j'}=\beta(\eta_{ij})\cdot g_{ij}\cdot \gamma_{j}$
 and that $\Phi'$ preserves identities follows from $\eta_{ii}=e$.
 That $\Phi'$ commutes with composition follows from
 $\eta_{ik}\cdot h_{ijk}=\gamma_{i}. h'_{\bar i'\bar j'\bar k'}\cdot \eta_{ij}\cdot g_{ij}.\eta_{jk}$
 (a corresponding equality of diagrams may be obtained from plugging
 together \eqref{eqn:fromCocyclesToBundles2} and \eqref{eqn:cobound2}).
 It is obvious from the definition that $\Phi$ commutes strictly
 with the previously defined action of $\cG$. Summarising, $\Phi$
 defines a strictly equivariant functor $\cP_{z}\to \cP_{z'}$ and,
 moreover, a morphism of principal 2-bundles.

 In the case that the coboundary is given by a common refinement
 $(V_{i})_{i\in I}$, which is properly finer than $(U_{i})_{i\in I}$ and
 $(U'_{i'})_{i'\in I'}$, we proceed as follows. Assume that
 $V_{i}\se U_{\bar i}$ and that $V_{i}\se U'_{\bar i'}$. Then we
 restrict $g_{\bar i\bar j}$ and $g'_{\bar i'\bar j'}$ to $V_{ij}$, as
 well as $h_{\bar i\bar j\bar k}$ and $h_{\bar i'\bar j'\bar k'}$ to
 $V_{\bar i'\bar j'\bar k'}$. This yields refined cocycles $\ol{z}$ and
 $\ol{z}'$ with the corresponding bundles $\cP_{\ol{z}}$,
 $\cP_{\ol{z}'}$ and the canonical inclusion
 $\Phi\from \cP_{\ol{z}}\to \cP_{z}$. Since $(V_{i})_{i\in I}$ is a
 refinement of itself, the previous construction yields morphisms
 $\cP_{\ol{z}}\to \cP_{\ol{z}'}$ and we set $\Phi'$ to be the
 composition of this morphism with the inclusion
 $\cP_{\ol{z}'}\to \cP_{z'}$.
\end{proof}
 
\begin{remark}\label{rem:morphismsAndOpenCovers}
 As the previous proposition suggests, $\check{H}(M,\cG)$ does not
 classify bundles in the classical sense, for bundle morphisms between
 different bundles may not be invertible. For instance, each cover
 $\mf{U}=(U_{i})_{i\in I}$ of $M$ gives rise to a cocycle with values in
 the trivial 2-group. For two different covers, all these cocycles are
 cohomologous, but one gets a morphism $\cP_{\mf{U}}\to \cP_{\mf{U}'}$
 (where $\cP_{\mf{U}}$, as constructed in Remark
 \ref{rem:fromCocyclesToBundles}, is simply the \v{C}ech-groupoid
 associated to $\mf{U}$) if and only if $\mf{U}$ is a refinement of
 $\mf{U}'$.

 However, if we start with a bundle $\cP$, extract a classifying cocycle
 $z(\cP)$ as in Proposition \ref{prop:fromBundlesToCocycles} (given by
 the choice of a cover $(U_{i})_{i\in I}$, trivialisations
 $\Phi_{i}, \ol{\Phi}_{i}$ and natural transformations $\ol{\tau}_{i}$)
 and reconstruct a bundle $\cP_{z(\cP)}$ as in Proposition
 \ref{prop:fromCocyclesToBundles}, then we always have a morphism
 $\coprod \ol{\Phi}_{i}\from \cP_{z(\cP)}\to \cP$. On objects, this
 morphism is given by $(i,x,g)\mapsto \ol{\Phi}_{i}(x,g)$ and on
 morphisms by the equivariant extension of
 \begin{equation}\label{eqn:morphismsAndopenCovers1}
  \left( (i,x,g)\xrightarrow{(i,j,x,e,e)}(j,x,g_{ij}^{-1})\right)\mapsto
  \left( \ol{\Phi}_{i}(x,e)\xrightarrow{\ol{\tau}_{j}(\ol{\Phi}_{i}(x,e))^{-1}}\ol{\Phi}_{j}(x,g_{ij}^{-1})\right)
 \end{equation}
 (with the notation from Remark \ref{rem:fromBundlesToCocycles}). That
 \eqref{eqn:morphismsAndopenCovers1} satisfies source-target matching
 follows from the definition of $g_{ij}$ (in the proof of Proposition
 \ref{prop:fromBundlesToCocycles}), implying
 \begin{equation*}
  \ol{\Phi}_{j}\circ \Phi_{j} \circ\ol{\Phi}_{i}(x,e)=\ol{\Phi}_{j} (x,g_{ij}^{-1})
 \end{equation*}
 on objects. That \eqref{eqn:morphismsAndopenCovers1} is also compatible
 with composition follows from the fact that $\coprod \ol{\Phi}_{i}$ is
 equivariant and from
 \begin{equation*}
  \Phi_{k}\Big(\ol{\tau}_{j}\big(\ol{\Phi}_{i}(x,e)\big)^{-1}\circ \ol{\tau}_{k}\big(\ol{\Phi}_{j}(x,g_{ij}^{-1})\big)^{-1}\Big)=\Phi_{k}\Big(\ol{\tau}_{k}\big(\ol{\Phi}_{i}(x,e)\big)^{-1}\cdot (h_{ijk},e)\Big),
 \end{equation*}
 which can be verified directly with the aid of
 \eqref{eqn:fromBundlesToCocycles_2}. Since
 \eqref{eqn:morphismsAndopenCovers1} is obviously smooth and equivariant
 by its definition, it defines a morphism $\cP_{z(\cP)}\to\cP$. Moreover,
 it follows from the fact that $\cG_{1}$ acts freely on $\cP_{1}$ that this
 functor is faithful.
\end{remark}

\begin{proposition}
 The smooth functor $\coprod\ol{\Phi}_{i}\from \cP_{z}\to \cP$, defined by
 \eqref{eqn:morphismsAndopenCovers1} is a weak equivalence (or strong Morita equivalence)
 of the underlying Lie groupoids.
\end{proposition}

\begin{proof}
 We have to show that
 \begin{itemize}
  \item [i)] the map
        \begin{equation*}
         \underbrace{\{((i,x,g),f)\in (\cP_{z})_{0}\times \cP_{1}:\Phi_{i}(x,g)=s(f)\}}_{:=Q}\xrightarrow{\ev} \cP_{0},\quad  ((i,x,g),f)\mapsto t(f)
        \end{equation*}
        admits local inverses and that
  \item [ii)] the diagram
        \begin{equation*}
         \xymatrix{
         \ar[d] \ar[r]
         \bigcup\,\,\{(i,j)\}\times U_{ij}\times H\times G	
         &
         \ar[d]
         \cP_{1}
         \\
         \ar[r]
         (\bigcup\,\, \{i\}\times U_{i}\times G)\times (\bigcup\,\, \{j\}\times U_{j}\times G)
         &
         \cP_{0}\times\cP_{0}
         }
        \end{equation*}
        is a pull-back.
 \end{itemize}
 To show i), we choose a local inverse of the target map in a local
 trivialisation and transform it to $\cP$. In fact, it is easily checked
 that for $((i,x,g),f)\in Q$ the map
 \begin{multline*}
  \pi^{-1}(U_{i})_{0}\ni y\mapsto \Big( \underbrace{(i,x(y),\beta(h_{f}^{-1}) g(y))}_{\in (\cP_{z})_{0}}, \\ \underbrace{\ol{\tau}_{i}(y)\circ \ol{\Phi}_{i}(x(y),(h_{f},\beta(h_{f}^{-1}) g(y)))\circ \ol{\tau}_{i}(\ol{\Phi}_{i}(x(y),\beta(h_{f}^{-1})g(y)))^{-1}}_{\in \cP_{1}} \Big),
 \end{multline*}
 where $x(y)$ and $g(y)$ denote the components of $\Phi_{i}(y)$ and
 $h_{f}$ is defined by $\Phi_{i}(f)=(x,h_{f},h_{f}^{-1}g(t(f)))$,
 defines a local (left) inverse for $\ev$, mapping $t(f)$ to
 $((i,x,g),f)$.

 In order to check ii), we verify the universal property directly. If
 $(i,x,g)$, $(j,y,k)$ and $f$ are given, such that $f$ is a morphism from
 $\ol{\Phi}_{i}(x,g)$ to $ \ol{\Phi}_{j}(y,k)$, then $f$ is a morphism
 in $\pi^{-1}(U_{i})$ and in $\pi^{-1}(U_{j})$, for both subcategories
 are full. Thus $\Phi_{i}(f)$ is a morphism from
 $(\Phi_{i}\circ \ol{\Phi}_{i})(x,g)=(x,g)$ to
 $(\Phi_{i}\circ \ol{\Phi}_{j})(y,k)=(y,g_{ji}(y)^{-1}k)$. From this it follows that
 $x=y$ and that $\Phi_{i}(f)=(x,(h_{f},g))$ with $\beta(h_{f})g=g_{ji}^{-1}(x)k$.
 Thus $(i,j,x,h_{iji}(x)^{-1}h_{f},g)$ is a morphism in $\cP_{z}$, which maps to
 $(i,x,g)$, $(j,x,k)$ and $f$ under the corresponding maps. This morphism
 is unique, because $\coprod\ol{\Phi}_{i}$ is faithful. Clearly, if
 $(i,x,g)$, $(j,x,k)$ and $f$ depend continuously on some parameter, then this morphism
 does also.
\end{proof}

\begin{definition}
 Two semi-strict principal $\cG$-2-bundles $\cP$
 and $\cP'$ over $M$ are said to be \emph{Morita equivalent} if there exists a
 third such bundle $\ol{\cP}$ and a diagram
 \begin{equation*}
  \vcenter{\xymatrix@=1em{
  & \ol{\cP}\ar[dl]_{\Phi} \ar[dr]^{\Phi'}\\
  \cP && \cP'
  }}
 \end{equation*}
 for semi-strict bundle morphisms $\Phi$ and $\Phi'$.
\end{definition}

\begin{lemma}
 Morita equivalence of bundles is in fact an equivalence relation.
\end{lemma}

\begin{proof}
 Suppose that we are given a diagram
 \begin{equation*}
  \vcenter{\xymatrix@=1em{
  & \ol{\cP}\ar[dl]_{\Phi} \ar[dr]^{\Phi'} & & \ol{\ol{\cP}}\ar[dl]_{\ol{\Phi}}\ar[dr]^{\bar{\Phi}'}\\
  \cP && \cP' && \cP''
  }}
 \end{equation*}
 implementing Morita equivalences between $\cP$ and $\cP'$, and between
 $\cP'$ and $\cP''$. Then Proposition \ref{prop:fromBundlesToCocycles}
 implies that $[z(\ol{\cP})]=[z(\ol{\ol{\cP}})]$ and thus there exists
 by Proposition \ref{prop:fromCocyclesToBundles} a bundle $\cQ$,
 together with morphisms $\cQ\to \cP_{z(\ol{\cP})}$ and
 $\cQ\to \cP_{z(\ol{\ol{\cP}})}$. 
 With the construction from Remark \ref{rem:morphismsAndOpenCovers} we can fill in the morphisms in the diagram
 \begin{equation*}
  \vcenter{\xymatrix@=1em{
	& \ar[d]\cP_{z(\ol{\cP})} &\ar[l]
\ar[dr]\ar[dl]
Q\ar[r]& \cP_{z(\ol{\ol{\cP}})}\ar[d] \\
  & \ol{\cP}\ar[dl]_{\Phi} & & \ol{\ol{\cP}}\ar[dr]^{\bar{\Phi}'}\\
  \cP &&  && \cP''
  }}
 \end{equation*}
 showing the claim.
\end{proof}

With this said, Proposition \ref{prop:fromBundlesToCocycles} and Proposition
\ref{prop:fromCocyclesToBundles} now imply the following classification theorem. 

\begin{theorem}
 If $\cG$ is a strict Lie 2-group and $M$ is a smooth manifold, then
 semi-strict principal $\cG$-2-bundles over $M$ are classified up to
 Morita equivalence by $\check{H}(M,\cG)$.
\end{theorem}

\begin{corollary}\label{cor:MoritaEquivalentLieGroupoids}
 If $\cP$ and $\cP'$ are Morita equivalent as principal 2-bundles, then the
 underlying Lie groupoids are also Morita equivalent.
\end{corollary}

We conclude this section with a couple of remarks on the classification result.

\begin{remark}
 Corollary \ref{cor:MoritaEquivalentLieGroupoids} shows in particular
 that the Morita equivalence class $[\cP]$ of a principal 2-bundle gives
 rise to a Morita equivalence class of the underlying Lie groupoid and
 thus determines a smooth stack. Moreover, the Lie 2-group $\cG$
 determines a group stack $[\cG]$. In fact, the Lie groupoid underlying
 $\cG$ can be given the structure of a stacky Lie group
 \cite{Blohmann08Stacky-Lie-groups} by turning the structure morphisms into
 bibundles as in \cite[Sect.\ 4.6]{Blohmann08Stacky-Lie-groups}, and
 this stacky Lie group gives a group stack, cf.\
 \cite{Blohmann08Stacky-Lie-groups}.

 Together with the morphism $[\pi]\from [\cP]\to [M]$, the right
 $[\cG]$-action on $[\cP]$ and the existence of local trivialisations
 give rise to something like a principal bundle in the 2-category of
 smooth stacks. One could have started our investigation with a
 rigourous definition of this concept and then pursuing a classification
 of those principal bundles in terms of non-abelian cohomology. This
 would also have lead to a classification in terms of non-abelian
 cohomology by very similar arguments. From this point of view it seems
 natural that non-abelian cohomology can classify principal 2-bundles
 only up to Morita equivalence. However, our approach is more direct and in
 more down-to-earth terms.
\end{remark}

\begin{remark}\label{rem:BandOfAGerbe}
 Let $(G,H,\alpha,\beta)$ be a smooth crossed module such that
 $\beta(H)$ is a normal split Lie subgroup of $G$ and let $\cG$ be the
 associated Lie 2-group. Then $G/\beta(H)$ carries a natural Lie group
 structure (cf.\ \cite[Def.\ 2.1]{neeb06nonAbelianExtensions}) and the
 projection map $G\to G/\beta(H)$ is smooth. If $\cP$ is a semi-strict
 principal $\cG$-2-bundle, then we obtain from this a principal
 $G/\beta(H)$-bundle $P$ by identifying isomorphic objects in $\cP$,
 i.e., we define $P$ to be $\Obj(\cP)/\sim$, where $p\sim p'$ if there
 exists a morphism between $p$ and $p'$. Then $P$ inherits naturally a
 $G/\beta(H)$-action, given by $[p].[g]:=[p.g]$ (where the dot between
 $p$ and $g$ refers to the $\cG$-2-space structure on $\cP$) and we
 endow $P$ with the quotient topology from $\Obj(\cP)$. Since $\cP\to M$
 (where $M$ is viewed as a category with only identity morphisms) maps
 isomorphic objects to the same element in $M$, this functor induces a
 map $P\to M$ (where $M$ is viewed as a space). If $(U_{i})_{i\in I}$ is
 an open cover such that there exist trivialisations
 $\Phi_{i}\from \left.\cP\right|_{U_{i}}\to U_{i}\times\cG$, then
 $\Phi_{i}$ induces an $G/\beta(H)$-equivariant bijective map
 $\left.P\right|_{U_{i}}\to U_{i}\times G/\beta(H)$ and we use this map
 to endow $P$ with a smooth structure. That this is in fact well-defined
 follows from the fact that the coordinate changes are then induced by
 the smooth maps $U_{ij}\to G/\beta(H)$, $x\mapsto [g_{ij}(x)]$, where
 $g_{ij}\from U_{ij}\to G$ is deduced from $\Phi_{i}$ and
 $\ol{\Phi}_{i}$ as in Remark \ref{rem:fromBundlesToCocycles}. This
 turns $P$ into a principal $G/\beta(H)$-bundle, which we call the
 \emph{band} of $\cP$.
\end{remark}

\begin{remark}\label{rem:groupoidExtensions}
 Another approach to assign differential geometric data to non-abelian
 \v{C}ech cohomology is to realise classes in $\check{H}(M,\cG)$ by
 Morita equivalence classes of Lie groupoid extensions, as outlined in
 \cite{Laurent-GengouxStienonXu09Non-abelian-differentiable-gerbes}. In
 particular, we recover
 \cite[3.14]{Laurent-GengouxStienonXu09Non-abelian-differentiable-gerbes}
 from the above classification by considering the crossed module
 $H\to G:=\Aut(H)$ (for a finite-dimensional $H$ with $\pi_{0}(H)$
 finite, say). For a non-abelian \v{C}ech cocycle $z$, Proposition
 \ref{prop:fromCocyclesToBundles} yields a 2-bundle $\cP_{z}$. Now $G$
 acts on the manifolds of objects and morphism of $\cP_{z}$ and since
 this action is obviously principal and all the structure maps of
 $\cP_{z}$ are compatible with the $G$-action, we have an induced Lie
 groupoid $\cP_{z}/G$, with objects $\coprod U_{i}$ and morphisms
 $\coprod U_{ij}\times H$. Moreover, the description of the composition
 in $\cP_{z}$ shows that $\cP_{z}/G$ is exactly the extension of
 groupoids from \cite[Prop.\
 3.14]{Laurent-GengouxStienonXu09Non-abelian-differentiable-gerbes}.

 However, $\cP_{z}/G$ is \emph{not} Morita equivalent to $\cP_{z}$. This
 can be seen for $M=\{*\}$, where $\cP_{z}$ is the action groupoid of
 $H$, acting via $\beta$ on $G$ and $\cP_{z}$ is the groupoid with one
 object and automorphism group $H$. Clearly, $\cP_{z}/G$ is transitive
 while $\cP_{z}$ is not.

 On the other hand, there is an extension of Lie groupoids, canonically
 associated to each principal 2-bundle, for an arbitrary
 finite-dimensional crossed module from now on. For this we note that
 the strong equivalence $\pi^{-1}(U_{i})\cong \ul{U_{i}}\times \cG$
 yields a weak equivalence \cite[Prop.\
 5.11]{MoerdijkMrcun03Introduction-to-foliations-and-Lie-groupoids} and
 we thus have
 \begin{multline*}
  \Mor(\pi^{-1}(U_{i}))\cong \{(p,p',(x,(h,g)))\in \cP_{0}\times\cP_{0}\times (U_{i}\times H\rtimes G):\\ \Phi_{i}(p)=(x,g),\,\Phi_{i}(p')=(x,\beta(h)\cdot g)\}
 \end{multline*}
 from the pull-back condition in the definition of weak equivalences.
 Moreover, the above diffeomorphism is in fact
 $(H\rtimes G)$-equivariant and we thus see that the action of
 $\ker(\beta)$ on $\Mor(\pi^{-1}(U_{i}))$ is principal. We thus have an
 associated extension
 \begin{equation}\label{eqn:LieGroupoidExtensionFromCocycle}
  \vcenter{\xymatrix@=1.5em{
  (\tx{identities in }\cP_{z})\cdot \ker(\beta)\ar[d]\ar[rr]&&  \cP_{1} \ar[rr]\ar@<-1ex>[d]\ar@<1ex>[d]&& \cP_{1}/\ker(\beta)\ar@<-1ex>[d]\ar@<1ex>[d]\\
  \cP_{0}\ar@{=}[rr]&&\cP_{0} \ar@{=}[rr]&&\cP_{0},
  }}
 \end{equation}
 of Lie groupoids. 

 Note also, that the construction of the band of a Lie groupoid extension from
 \cite{Laurent-GengouxStienonXu09Non-abelian-differentiable-gerbes} differs
 from the construction in Remark \ref{rem:BandOfAGerbe}, for the band there
 is a principal bundle over the space of objects of the considered Lie groupoid,
 while the band that we construct is a principal bundle over the quotient
 $\cP_{0}/\cP_{1}$, if it exists as a manifold.

 It would be interesting to understand the exact
 correspondence between our approach and
 \cite{Laurent-GengouxStienonXu09Non-abelian-differentiable-gerbes} in
 more detail.
\end{remark}

\begin{remark}
 If $\beta (H)$ is a split Lie subgroup, so that $K:=G/\beta (H)$ is
 again a smooth Lie group, then each non-abelian cocycle determines
 a smooth $K$-valued 1-cocycle $k_{ij}\from U_{ij}\to K$ and because
 cohomologous cocycles are mapped to cohomologous 1-cocycles, we thus
 get a map
 \[
  Q\from \check{H}(M,\cG)\to\check{H}^{1}(M,K)
 \]
 (realised on bundles by the preceding construction). This map is
 surjective, for each $U_{ij}$ is contractible, and thus each map
 $U_{ij}\to K=G/\beta (H)$ has a lift to $G$. The fibres of this map
 then classify semi-strict principal $\cG$-2-bundles with a fixed
 underlying band. In particular, if $H$ is abelian, the fibre of the
 trivial band (i.e., all $g_{ij}$ take values in $\beta (H)$) is
 isomorphic to $H^{2}(M,H)$.
\end{remark}

\begin{remark}
 One can also define a topological version $\check{H}_{\op{top}}(M,\cG)$
 of $\check{H}(M,\cG)$, where all occurring maps $g_{ij}$, $h_{ijk}$,
 $\gamma_{i}$ and $\eta_{ij}$ are required to be continuous rather than
 smooth. The same classification goes through along the same lines for
 topological $\cG$-2-bundles (for $\cG$ a topological 2-group) over an
 arbitrary topological space $M$. If $M$ is paracompact, then
 $\check{H}_{\op{top}}(M,\cG)$ stands in bijection with the set of
 homotopy classes $[M,B|N\cG|]$ and consequently with
 $\check{H}_{\op{top}}^{1}(M,|N\cG|)$, where $|N\cG|$ is a topological
 group, associated to $\cG$ (the geometric realisation of the nerve of
 the category $\cG$). This has been shown in
 \cite{BaezStevenson08The-Classifying-Space-of-a-Topological-2-Group},
 cf.\ also 
 \cite{Jurco05Crossed-Module-Bundle-Gerbes;-Classification-String-Group-and-Differential-Geometry}
 (note that $\cG$ is always well-pointed in our case, for we are only
 dealing with Lie groups). In particular, $\check{H}_{\op{top}}(M,\cG)$
 is trivial if $M$ is paracompact and contractible. This shows that for
 paracompact finite-dimensional $M$, one can always assume that bundles
 are trivialised over a fixed good cover and one does not run into the
 problems described in Remark \ref{rem:morphismsAndOpenCovers}. A
 similar approach as in
 \cite{MullerWockel07Equivalences-of-Smooth-and-Continuous-Principal%
 -Bundles-with-Infinite-Dimensional-Structure-Group} should yield the
 same for $\check{H}(M,\cG)$.
\end{remark}

\begin{remark}\label{rem:classificationForCentralExtensions}
 Assume that $\beta$ is surjective (i.e., assume that $\beta$ is a
 central extension) and set $\check{\cH}:=(\{*\},\ker(\beta))$. Then $H(M,\cG)$
 is isomorphic (as a set) to $\check{H}(M,\cH)$ for paracompact and
 finite-dimensional $M$. In fact, if $(g_{ij},h_{ijk})$ is a
 non-abelian cocycle, then we define a cohomologous cocycle as
 follows. First, we assume w.l.o.g. that each $U_{ij}$ is contractible,
 so that $g_{ij}$ lifts to $\eta_{ij}\from U_{ij}\to H$ (assuming
 $\eta_{ii}\equiv e_{H}$), and we set $\gamma_{i}$ to be constantly
 $e_{G}$. Then \eqref{eqn:cobound1} and \eqref{eqn:cobound2} define a
 cohomologous cocycle $(g'_{ij},h'_{ijk})$ and from \eqref{eqn:cocyc1}
 it follows that $g'_{ij}$ is also constantly $e_{G}$ and thus
 $h'_{ijk}$ takes values in $\ker(\beta)$. Thus the canonical map
 $\check{H}(M,\cH)\to H(M,\cG)$ is surjective and the injectivity
 follows similarly. Consequently, principal $\cG$ 2-bundles are
 classified (up to Morita equivalence) by $H^{2}(M,\ker(\beta ))$ if
 $\beta$ is surjective.
\end{remark}

\begin{remark}
 There is also a way of understanding the construction in Proposition
 \ref{prop:fromCocyclesToBundles}, given by a construction of 2-bundles
 as quotients of equivalence 2-relations as in \cite[Prop.\
 22]{bartels06HigherGaugeTheoryI:2-Bundles}.

 Let $\mathcal{U}=(U_{i})_{i\in I}$ be an open cover of $M$.  We set
 $Y:=\coprod U_{i}$ and $\pi_{1}\from Y\to M$, $(i,x)\mapsto x$. Then
 $Y^{[n]}:=Y\times_{M}\dots \times_{M}Y$ (n-fold fibre product) is the
 disjoint union of $n$-fold intersections of the $U_{i}$ and we denote by
 $\mathcal{U}^{[n]}$ the corresponding category with only identity
 morphisms. Moreover, we have canonical projections
 $\pi_{n_{1}\dots n_{k}}\from Y^{[n]}\to Y^{[n-k]}$ for $k<n$, which we
 identify with the corresponding functors
 $\pi_{n_{1}\dots n_{k}}\from\mathcal{U}^{[n]}\to \mathcal{U}^{[n-k]}$.

 A non-abelian cocycle $c=(g_{ij},h_{ijk})$, with underlying open
 cover $\mathcal{U}$, defines what is called a 2-transition in
 \cite[Sect.\ 2.5.1]{bartels06HigherGaugeTheoryI:2-Bundles}. The functor
 $g\from \mathcal{U}^{[2]}\to \cG$ (called 2-map in
 \cite{bartels06HigherGaugeTheoryI:2-Bundles}) is given by the smooth
 map $g\from Y^{[2]}\to G$, $((i,j),x)\mapsto g_{ij}(x)$, and the
 natural isomorphism
 \[
  \gamma\from \mu \circ (g\times g)\circ (\pi_{01}\times
  \pi_{12})\Rightarrow g\circ \pi_{02}
 \]
 is then given by
 $Y^{[3]}\ni ((i,j,k),x)\mapsto (h_{ijk}(x),g_{ij}(x)\cdot g_{jk}(x))\in H\rtimes G$.
 This 2-transition is semi-strict in the sense of
 \cite{bartels06HigherGaugeTheoryI:2-Bundles} (i.e., the natural
 $\mu \circ g\circ \iota \Rightarrow \idObject$ for
 $\iota \from \mathcal{U}\hookrightarrow \mathcal{U}[2]$,
 $(i,x)\mapsto ((i,i),x)$ is the identity), for $g_{ii}\equiv e_{G}$ in
 our setting. Note that $\gamma$ is a natural isomorphism because of
 condition \eqref{eqn:cocyc1} and the coherence, required in
 \cite[Sect.\ 2.5.1]{bartels06HigherGaugeTheoryI:2-Bundles} is condition
 \eqref{eqn:cocyc2}).

 In \cite[Prop.\ 22]{bartels06HigherGaugeTheoryI:2-Bundles}, the bundle
 is constructed from the 2-transition $(g,\gamma )$ by taking the
 quotient of the category $\mathcal{U}^{[2]}\times \cG$ by an
 equivalence 2-relation, determined by $(g,\gamma )$. This equivalence
 2-relation is a categorified version of an equivalence relation,
 expressed purely in arrow-theoretical terms (cf.\ \cite[1.1.4 and
 2.1.4]{bartels06HigherGaugeTheoryI:2-Bundles}).

 This equivalence 2-relation is determined by two functors
 $\mathcal{U}^{[2]}\times \cG\to\mathcal{U}\times\cG$, one given by
 $\pi_{1}\times \id_{\cG}$ and the other one by
 \[
  (\id_{\mathcal{U}}\times \rho )\circ (\id_{\mathcal{U}}\times g\times
  \id_{\cG})\circ \iota\times \id_{\cG}.
 \]
 One readily checks that these two functors are what is called jointly
 2-monic in \cite{bartels06HigherGaugeTheoryI:2-Bundles}, for natural
 equivalences are basically given by $H$-valued mappings, allowing
 lifts of natural equivalences to be constructed directly. The
 2-reflexivity map is given by $\iota\times \id_{\cG}$ (and identities
 as natural isomorphisms, for our 2-transition is semi-strict). The
 2-kernel pair of $\pi_{2}\times\id_{\cG},\pi_{1}\times \id_{\cG}$ is
 simply
 \[
  \begin{CD}
   \mathcal{U}^{[3]}\times \cG @>{\pi_{23}\times \id_{\cG}}>>
   \mathcal{U}^{[2]}\times \cG\\
   @V{\pi_{12}\times \id_{\cG}}VV @VV{\pi_{1}\times \id_{\cG}}V\\
   \mathcal{U}^{[2]}\times \cG @>{\pi_{2}\times \id_{\cG}}>>
   \mathcal{U}\times\cG.
  \end{CD}
 \]
 The Euclideanness functor is given by $\pi_{13}\times \id_{\cG}$ and
 the first equivalences in the Euclideanness condition is trivial and
 the second one is given by $\gamma$ (we choose the 2-kernel pair to be
 defined by $\pi_{1}$ and $\pi_{0}$ so that it fits with the usual
 notion of an equivalence relation).

 It can be checked that the category $\cP_{c}$ is a quotient of this
 equivalence relation by the inclusion
 $\mathcal{U}\times \cG \hookrightarrow \cP_{c}$ (cf.\ \cite[Sect.\
 2.1.4]{bartels06HigherGaugeTheoryI:2-Bundles}) and thus realises the
 bundle constructed in \cite[Prop.\
 22]{bartels06HigherGaugeTheoryI:2-Bundles}. We leave the details as an
 exercise. From this construction one sees immediately that the quotient
 exists in the category of smooth manifolds.
\end{remark}

\section{Gauge 2-groups}

In the classical setup, a gauge transformation of a principal bundle is
a bundle self-equivalence and all gauge transformations form a group
under composition. Likewise, in the categorified case the vertical
self-equivalences form a category (as functors and natural
transformations) which is in fact a weak 2-group with respect to the natural
compositions.

We will show that this weak 2-group is in fact equivalent to a naturally given
strict 2-group. Moreover, we show that under some mild conditions, this strict
2-group carries naturally the structure of a strict \emph{Lie} 2-group. As in the
previous section, the fact that we only consider strict actions shall be the crucial
point to make the ideas work.

Unless stated otherwise, we assume throughout this section that
$\cG$ is a strict Lie 2-group arising from the smooth crossed module
$(\alpha ,\beta ,G,H)$ and that $\cP$ is a semi-strict principal
2-bundle over the smooth manifold $M$. We will identify $M$ with the smooth
2-space it determines by adding only identity morphisms.

\begin{remark}
 We consider the category $\Aut(\cP)^{\cG}$, whose objects are morphisms
 $F\from \cP\to\cP$ of principal $\cG$-2-bundles and whose morphisms are
 2-morphisms $\alpha\from F\Rightarrow G$ (cf.\ Definition
 \ref{def:principal-2-bundle}). This is a weak 2-group with respect to
 composition of functors and natural equivalences (cf.\ \cite[Ex.\
 34]{baezLauda04Higher-DimensionalAlgebraV:2-Groups}). We call this weak
 2-group the \emph{gauge 2-group} of $\cP$.
\end{remark}

We shall make this weak 2-group more accessible by showing that it is
equivalent to a strict 2-group. Initially, we start with the most simple
case.

\begin{proposition}
 The category $\Aut(\cG)^{\cG}$ of strictly equivariant endofunctors of
 $\cG$ is equivalent to $\cG$.
\end{proposition}

\begin{proof}
 Each $F\in\Aut(\cG)^{\cG}$ is given by a functor $F\from \cG\to\cG$
 satisfying $F(g\cdot g')=F(g)\cdot g'$ on objects and
 $F((h,g)\cdot (h',g'))=F((h,g))\cdot (h',g')$ on morphisms. From this
 it follows that $F(h,g)=(k_{1}\cdot k_{2}.h,k_{2}\cdot g)$, where
 $F((e,e))=(k_{1},k_{2})$, and the compatibility with the structure maps
 yields $k_{1}=e$. Likewise, a natural equivalence between such functors
 is uniquely given by its value at $e_{G}$, which is an element of $H$.
\end{proof}

Recall that for a category $\cC$, we denote by
$\Delta\from \cC\to\cC\times\cC$ the diagonal embedding and $\Delta_{0}$
denotes its map on objects.

\begin{lemma}\label{lem:mappingGroupIsA2-Group}
 Let $\cM$ be an arbitrary smooth
 2-space. Then the category $\cC^{\infty}(\cM,\cG)$ of smooth functors
 from $\cM$ to $\cG$ is a 2-group with respect to the monoidal functor
 $\mu_{*}$, given on objects by
 $\mu_{*}(F,G):= \mu \circ (F\times G)\circ \Delta$ and on morphisms by
 $\mu_{*}(\alpha,\beta ):= \mu_{1}\circ (\alpha \times \beta)\circ \Delta_{0}$.
 Moreover, if $\cG$ is strict, then $\cC^{\infty}(\cM,\cG)$ is so.
\end{lemma}

\begin{proof}
 To check that $\mu_{*}(\alpha ,\beta)$ is a natural transformation from
 $\mu_{*}(F,F')$ to $\mu_{*}(G,G')$, one computes that it coincides with
 the horizontal composition of the natural transformations
 \[
  (\id\from \mu \Rightarrow \mu)\circ ((\alpha ,\beta)\from
  (F,F')\Rightarrow(G,G'))\circ (\id\from\Delta\Rightarrow \Delta).
 \]
 The rest is obvious.
\end{proof}

\begin{remark}\label{rem:push-forwardCrossedModule}
 One can easily read off from Lemma \ref{lem:mappingGroupIsA2-Group} the
 crossed module $(\alpha_{*},\beta_{*},G_{*},H_{*})$, that yields
 $\cC^{\infty}(\cM,\cG)$ as strict 2-group (cf.\ Remark
 \ref{rem:2-groupsFromCrossedModules}). The objects of
 $\cC^{\infty}(\cM,\cG)$ form a set $\cMor(\cM,\cG)$ (as a subset of
 $C^{\infty}(\cM_{0},\cG_{0})\times C^{\infty}(\cM_{1},\cG_{1})$) and
 $\mu_{*}$ defines a group multiplication on this set. Thus we set
 $G_{*}=\cMor (\cM,\cG)$. Moreover, it is easily checked that
 $\cMor (\cM,\cG)$ is a subgroup of
 $C^{\infty}(\cM_{0},\cG_{0})\times C^{\infty}(\cM_{1},\cG_{1})$.
 Likewise, the morphisms in $\cC^{\infty}(\cM,\cG)$ form a set $2$-%
 $\cMor(\cM,\cG)$ and $\mu_{*}$ defines a group multiplication on this
 set. Again, one can interpret $2$-$\cMor(\cM,\cG)$ as a subgroup
 \[
  2\text{-}\cMor(\cM,\cG)\leq \cMor(\cM,\cG)\times
  C^{\infty}(\cM_{0},\cG_{1})\times \cMor(\cM,\cG),
 \]
 with $((F,\alpha ,G)\in 2$-
 $\cMor(\cM,\cG)):\Leftrightarrow(\alpha \from F\Rightarrow G)$. We set
 $H_{*}$ to be the kernel of the source map as a subgroup of
 $C^{\infty}(\cM_{0},\cG_{1})\times \cMor(\cM,\cG)$.  Then the
 homomorphism $\beta_{*} \from H_{*}\to G_{*}$ is the projection to the
 second component and the action $\alpha_{*}$ of $G_{*}$ on $H_{*}$ is
 the conjugation action on the second component.

 If $\cM=(M,M)$ has only identity morphisms, then
 $\cMor(\cM,\cG)\cong C^{\infty}(M,\cG_{0})$ and
 $2$-$\cMor (\cM,\cG)\cong C^{\infty}(M,\cG_{1})$. From this
 it follows that $\cC^{\infty}(\cM,\cG)$ is associated to
 the push-forward crossed module
 $(\alpha_{*},\beta_{*},C^{\infty}(M,G),C^{\infty}(M,H))$, where
 $\alpha_{*}$ and $\beta_{*}$ are the point-wise applications of
 $\alpha$ and $\beta$.
\end{remark}

The following proposition can be understood as an instance of the fact
that in the classical case, a bundle endomorphism (covering the identity
on the base) of a principal bundle is automatically invertible, and thus
bundle endomorphisms form a group. This can best be verified by viewing
bundle maps as smooth group-valued maps on the total space. However,
note that morphisms between distinct principal bundles need not be
invertible (cf.\ Remark \ref{rem:morphismsAndOpenCovers}).

\begin{proposition}\label{prop:gauge2-GroupIsMapping2-Group}
 The weak 2-group $\Aut(\cP)^{\cG}$ of self-equivalences of $\cP$ is
 equivalent, as a weak 2-group, to $\cC^{\infty}(\cP ,\cG_{\Ad})^{\cG}$,
 the strict 2-group of morphisms of $\cG$-2-spaces, where $\cG_{\Ad}$
 denotes $\cG$ with the conjugation action from the right.
\end{proposition}

\begin{proof}
 The existence of strictly equivariant local trivialisations imply that
 $\cP_{x}:=\pi^{-1}(x)$ is equivalent to $\cG$. Then the usual reasoning
 gives
 $\Aut(\cP_{x},\cP_{x})^{\cG}\cong \cFun(\cP_{x},\cG_{\Ad})^{\cG}$. Since
 self-equivalences preserve the
 subcategories $\cP_{x}$, each object in $\Aut(\cP)^{\cG}$ is thus
 given by a strictly equivariant functor
 $\gamma_{F}\from \cP\to \cG_{\Ad}$. That this functor is in fact smooth
 can be seen in local coordinates. Likewise, each smooth 2-morphism
 $\alpha \from F\Rightarrow G$ between morphisms $F$ and $G$ of
 $\cG$-2-spaces is given by a smooth equivariant map
 $\eta_{\alpha } \from \cP_{0}\to H\rtimes G$.

 It is readily checked that $F\mapsto \gamma_{F}$ and
 $\alpha\mapsto \eta_{\alpha}$ defines a monoidal functor from
 $\Aut(\cP)^{\cG}$ to $\cC^{\infty}(\cP,\cG)^{\cG}$.  The inverse
 functor from $\cC^{\infty}(\cP,\cG)^{\cG}$ to $\Aut(\cP)^{\cG}$ is
 obviously given by $\gamma \mapsto F_{\gamma }$ on objects and
 $\eta \mapsto \alpha_{\eta}$ on morphisms, where
 $F_{\gamma}=\rho \circ (\id_{\cP}\times \gamma)\circ \Delta_{\cP}$ and
 $\alpha_{\eta}(p_{0})=\id_{p_{0}}\cdot \eta (p_{0})$.
\end{proof}

\begin{remark}
 For a semi-strict principal 2-bundle $\cP_{c}$, given by a
 non-abelian cocycle $c=(h_{ijk},g_{ij})$, we can also interpret the
 equivalence $\Aut(\cP_{c})^{\cG}\cong \cC^{\infty}(\cP_{c},\cG)^{\cG}$
 as follows. As we have seen in the proof of Proposition
 \ref{prop:fromBundlesToCocycles}, each self-equivalence of $\cP_{c}$
 gives rise to an equivalence of $c$, given by smooth maps
 $\gamma_{i}\from U_{i}\to G$ and $\eta_{ij}\from U_{ij}\to H$, obeying
 \eqref{eqn:cobound1}-\eqref{eqn:cobound2} and normalisation. This
 defines a smooth functor $\cP_{c}\to\cG$, given on objects by
 $(i,x,g)\mapsto g\cdot \gamma_{i}(x)$ and on morphisms by
 $((i,j),x,(h,g))\mapsto (\gamma_{j}(x)^{-1}.(h\cdot \eta_{ij}(x)),g\cdot \gamma_{j}(x))$.
\end{remark}

We now turn to the smoothness conditions on $\Aut(\cP)^{\cG}$. We will
endow all spaces of continuous maps with the smooth $C^{\infty}$
topology, i.e., if $M$ and $N$ are smooth manifolds, then we endow
$C^{\infty}(M,N)$ with the initial topology with respect to
\begin{equation*}
 C^{\infty}(M,N)\to \prod_{k=0}^{\infty}C(T^{k}M,T^{k}N),\quad f\mapsto (T^{k}f)_{k\in \N_{0}}
\end{equation*}
(where $C(T^{k}M,T^{k}N)$ is equipped with the compact-open topology).
This is the topology on spaces of smooth functions used in \cite{autP}
and \cite{centrExt}, whose results we shall cite in the sequel in order
to establish Lie 2-group structures on gauge 2-groups.

\begin{proposition}
 If $M$ is compact, then $\cC^{\infty}(M,\cG)$ is a Lie 2-group, which
 is associated to the smooth crossed module
 $(\alpha_{*},\beta_{*},C^{\infty}(M,G),C^{\infty}(M,H))$.
\end{proposition}

\begin{proof}
 We have already seen in Remark \ref{rem:push-forwardCrossedModule} that
 $\cC^{\infty}(M,\cG)$ is associated to the push-forward crossed module
 $(\alpha_{*},\beta_{*},C^{\infty}(M,G),C^{\infty}(M,H))$. This is
 actually a smooth crossed module, the only non-trivial thing to check
 is the smoothness of the action of $C^{\infty}(M,G)$ on
 $C^{\infty}(M,H)$. But this follows from the smoothness of
 parameter-dependent push-forward maps (cf.\ \cite[Prop.\
 3.10]{gloeckner02a} and \cite[Prop.\ 28]{smoothExt}) and the fact that
 automorphic actions need only be smooth on unit neighbourhoods in order
 to be globally smooth.
\end{proof}

Before coming to the main result of this section, we have to provide
some Lie theory for strict Lie 2-groups by hand.

\begin{remark}
 We briefly recall strict Lie 2-algebras
 \cite{baezCrans04HigherDimensionalAlgebraVILie2-Algebras}. The
 definition is analogous to that of a strict Lie 2-group as a category in
 Lie groups. First, a 2-vector space is a category, in which all spaces
 are vector spaces and all structure maps
 are linear. A strict Lie 2-algebra is
 then a 2-vector space $\fG$, together with a functor
 $[\cdot ,\cdot ]\from \fG\times\fG\to\fG$, which is required to be
 linear and skew symmetric on objects and morphisms and which satisfies
 the Jacobi identity
 \[
  [x,[y,z]]=[[x,y],z]+[y,[x,z]]
 \]
 on objects and morphisms.

 Coming from strict Lie 2-groups, there is a natural way to associate a
 strict Lie 2-algebra to a strict Lie 2-group by applying the Lie
 functor $G\mapsto T_{e}(G)$, $f\mapsto Tf(e)$. This works, because this
 functor preserves pull-backs and thus all categorical structures (cf.\
 \cite[Prop.\ 5.6]{baezCrans04HigherDimensionalAlgebraVILie2-Algebras}).
 If $\cG$ is a Lie 2-group, then we denote by $L(\cG)$ the strict Lie
 2-algebra one obtains in this way.

 We have the same interplay between crossed modules of Lie algebras and
 strict Lie 2-algebras as in the case of strict Lie 2-groups. A \textit{crossed
 module} (of Lie algebras) consists of two Lie algebras $\fg,\fh$, and
 action $\dot{\alpha}\from \fg\to\der(\fh)$ and a homomorphism
 $\dot{\beta}\from\fh\to\fg$ satisfying
 $\dot{\beta}(\dot{\alpha}(x).y)=[x,\dot{\beta}(y)]$ and
 $\dot{\alpha}(\dot{\beta}(x)).y=[x,y]$. To such a crossed module one
 can associate the Lie 2-algebra $(\fg,\fh\rtimes \fg)$ with $s(x,y)=y$,
 $t(x,y)=\dot{\beta}(x)+y$, $(z,\dot{\beta}(x)+y)\circ (x,y)=(z+x,y)$
 and $[\cdot ,\cdot]$ given by the Lie-bracket on $\fg$ and
 $\fh\rtimes\fg$. Moreover, one checks readily that if $\cG$ is
 associated to $(\alpha ,\beta ,G,H)$, then $L(\cG)$ is associated to
 the derived crossed module $(\dot{\alpha},\dot{\beta},\fh,\fg)$.
\end{remark}

\begin{theorem}\label{thm:LiesThirdTheorem}
 If $\fG$ is a strict Lie 2-algebra with finite-dimensional object- and morphism space,
 then there
 exists a strict Lie 2-group $\cG$ such that $L(\cG)$ is isomorphic to
 $\fG$.
\end{theorem}

\begin{proof}
 There is a functor from Lie algebras to simply connected Lie groups,
 which is adjoint to the Lie functor. This functor also preserves
 pull-backs and applied to the spaces of objects and morphisms and to
 the structure maps of a strict Lie 2-algebra produces a strict Lie
 2-group.
\end{proof}

\begin{remark}
 If $\cG$ is a strict Lie 2-group with strict Lie 2-algebra $\fG$, then
 we also have a strict 2-action $\cAd\from \fG\times \cG\to\fG$ of $\cG$
 on $\fG$. This is given on objects and morphisms by
 $(x,g)\mapsto \Ad(g^{-1}).x$, where $\Ad$ is the ordinary adjoint
 action. That this defines a functor follows from
 $\Ad(\varphi (g)).\dot{\varphi}(x)=\dot{\varphi}(\Ad (g).x)$ for each
 homomorphism $\varphi$ of Lie groups. We denote the corresponding
 $\cG$-2-spaces by $\cG_{\Ad}$ and $\fG_{\op{ad}}$.
\end{remark}

The Lie 2-algebra which is of particular interest in this section is the
following.

\begin{proposition}
 If $\cM$ is a strict $\cG$-2-space, then
 $\cC^{\infty}(\cM,L(\cG)_{\op{ad}})^{\cG}$, the category of morphisms of
 $\cG$-2-spaces from $\cM$ to $L(\cG)_{\op{ad}}$
 is a strict Lie 2-algebra. The functor $[\cdot ,\cdot]$ is given by the
 point-wise application of the functor in $L(\cG)$ as in Lemma
 \ref{lem:mappingGroupIsA2-Group}.
\end{proposition}

\begin{proof}
 We set $\fG:=L(\cG)$, $\fK:=\cC^{\infty}(\cM,\fG)^{\cG}$ and identify
 $\fK_{0}$ with a subset of the locally convex Lie algebra
 $\fl:=C^{\infty}(\cM_{0},\fG_{0})\times C^{\infty}(\cM_{1},\fG_{1})$.
 The requirement on $(\xi ,\nu)$ to define a functor may be expressed in
 terms of point evaluations and linear maps, for instance, the
 compatibility with the source map is
 \[
  s_{\fG}(\nu(x)) =\xi(s_{\cM}(x))\fa x\in \cM_{1}.
 \]
 The same argument applies for the requirement on a functor to be
 $\cG$-equivariant, and thus $\fK_{0}$ is a closed subalgebra in $\fl$.

 In the same way, we may view $\fK_{1}$ as a closed subalgebra of
 $ \fK_{0}\times C^{\infty}(\cM_{0},\fG_{1})\times \fK_{0}$,  with
 $(F,\alpha ,G)\in \Mor(\cC^{\infty}(\cM,\fG))$ if and only if
 $\alpha \from F\Rightarrow G$ is a smooth natural equivalence. The
 structure maps are given by projections, embeddings and push-forwards
 by continuous linear mappings and thus all continuous algebra
 morphisms.
\end{proof}

The crucial tool in the description of the Lie group structure on
$\cC^{\infty}(\cP,\cG_{\Ad})^{\cG}$ shall be the exponential functions on $G$
and $H\rtimes G$ in the case that it provides charts for the Lie group
structures (i.e., if $G$ and $H\rtimes G$ are \textit{locally
exponential}, cf.\ Appendix
\ref{app:differentialCalculusOnSpacesOfMappings}).

\begin{lemma}
 If $\cG$ is a strict Lie 2-group, such that its group of objects and
 morphisms possess an exponential function, then these functions define
 a smooth functor
 \[
  \cExp\from L(\cG)\to\cG.
 \]
\end{lemma}

\begin{proof}
 For each homomorphism $\varphi \from G_{1}\to G_{1}$ between Lie groups
 with exponential function, the diagram
 \begin{equation}\label{eqn:exponentialMapCommutesWithMorphisms}
  \begin{CD}
  G_{1}@>{\varphi }>>G_{2}\\
  @A{\exp_{1}}AA @A{\exp_{2}}AA\\
  \fg_{1}@>{L(\varphi)}>>\fg_{2}
  \end{CD}
 \end{equation}
 commutes. Since all requirements on $\cExp$ to define a functor can be
 phrased in such diagrams, the assertions follows.
\end{proof}

\begin{remark}\label{rem:equivarianceOf2-Morphisms}
 Let $\cM$ be a strict $\cG$-2-space. If the 2-morphism
 $\alpha\from F\Rightarrow G$ between the morphisms $F,G$ of
 $\cG$-2-spaces is viewed as a map $\alpha\from \cM_{0}\to H\rtimes G$,
 then $\alpha=(\ol{\alpha},F_{0})$ for some
 $\ol{\alpha}\from \cM_{0}\to H$ and $G$ satisfies
 satisfies $G_{0}=(\beta \circ \ol{\alpha})\cdot F_{0}$ and
 \[
  G_{1}=((\ol{\alpha}\circ t_{\cM})\cdot \ol{F_{1}}\cdot (\alpha \circ
  s_{\cM})^{-1},G_{0}\circ s_{\cM}).
 \]
 Thus $G$ is uniquely determined by $F$ and $\ol{\alpha}$. If $F$ and
 $\ol{\alpha}$ are is strictly equivariant, then so is $G$ since
 $\ol{\alpha}\in C^{\infty}(\cM_{0},H)^{G}$ by definition.
\end{remark}

For classical principal bundles, the compactness of the base manifold
and the local exponentiality of the structure group ensure the existence
of Lie group structures on gauge transformation groups (cf.\
\cite{autP}). We shall follow similar ideas here and call a
strict Lie 2-group \textit{locally
exponential} if its Lie groups of objects and morphisms are so.

\begin{theorem}\label{thm:gauge2-GroupIsLie2-Group}
 Assume that $M$ is compact, that $\cG$ is locally exponential, and that
 the actions of $\cG_{1}$ on $\cP_{1}$ and of $\cG_{0}$ on $\cP_{0}$ are
 principal. Then $\cC^{\infty}(\cP,\cG_{\Ad})^{\cG}$ is a locally exponential
 strict Lie 2-group with strict Lie 2-algebra $\cC^{\infty}(\cP,L(\cG)_{\op{ad}})^{\cG}$.
\end{theorem}

\begin{proof}
 The proof works similar as in the case for classical principal bundles
 in \cite{autP}. We set $\cK :=\cC^{\infty}(\cP,\cG)^{\cG}$,
 $\fK:=\cC^{\infty}(\cP,L(\cG))^{\cG}$ and denote by $\cK_{0}$,
 $\fK_{0}$ and $\cK_{1}$, $\fK_{1}$ the corresponding spaces of objects
 and morphisms.  Then we have
 \begin{multline*}
  \cK_{0}=\{(\gamma ,\eta)\in C^{\infty}(\cP_{0},G)^{G}\times
  C^{\infty}(\cP_{1},H\rtimes G)^{H\rtimes G}: 
  \eta \circ i_{\cP}=i_{\cG}\circ \gamma \\
  \gamma \circ s_{\cP}=s_{\cG}\circ \eta, \,\,
  \gamma \circ t_{\cP}=t_{\cG}\circ \eta, \,\, 
  \eta \circ c_{\cP}= c_{\cG}\circ (\eta _{\,s\!\!}\times_{t\!} \eta) 
  \}.
 \end{multline*}
 Since the conditions on $(\gamma ,\eta)$ to be in $\cK_{0}$ can all be
 phrased in terms of evaluation maps on $\cP_{0}$ and $\cP_{1}$, it
 follows that $\cK_{0}$ is a closed subgroup of
 $L:=C^{\infty}(\cP_{0},G)^{G}\times C^{\infty}(\cP_{1},H\rtimes G)^{H\rtimes G}$,
 endowed with the $C^{\infty}$-topology. Similarly, we obtain $\fK_{0}$
 as a closed subalgebra of
 $\fl:=C^{\infty}(\cP_{0},\fg)^{G}\times C^{\infty}(\cP_{1},\fh\rtimes \fg)^{H\rtimes G}$.

 Now $L$ is a locally convex, locally exponential Lie group, modelled on
 $\fl$ (cf.\ \cite[Thm.\ 1.11]{autP}), because the actions of $G$ on
 $\cP_{0}$ and of $H\rtimes G$ on $\cP_{1}$ are free and locally trivial
 (cf.\ Remark \ref{rem:groupoidExtensions}). The exponential function
 for this Lie group is then given by
 \[
  \fl\ni(\xi ,\nu )\mapsto (\exp_{G}\circ \xi ,\exp_{H\rtimes G}\circ
  \nu )\in L,
 \]
 which restricts to a diffeomorphism on some zero neighbourhood in
 $\fl$. Since this is the same as the composition of the exponential
 functor $\cExp$ with $(\xi ,\nu)$, this exponential function restricts
 to a map from $\fK_{0}$ to $\cK_{0}$. It follows from the construction
 of the Lie group structure on $L$ that this map restricts to a
 bijective map of an open zero neighbourhood in $\fK_{0}$ to an open
 unit neighbourhood in $\cK_{0}$. The compatibility with the structure
 maps of $\cP$, $L(\cG)$ and $\cG$ can be checked by repeated use of
 \eqref{eqn:exponentialMapCommutesWithMorphisms} in local coordinates,
 for all structure maps of $L(\cG)$ and $\cG$ are (by the construction
 of $L(\cG)$) given by morphisms of Lie algebras and Lie groups,
 commuting with the respectively exponential functions. For instance
 \[
  \exp_{G}\circ \xi \circ s_{\cP}=s_{\cG}\circ \exp_{H\rtimes G}\circ
  \nu \quad\Leftrightarrow\quad \xi \circ s_{\cP}=s_{L(\cG)}\circ \nu.
 \]
 if $\xi$ and $\nu$ have representatives in local coordinates, which
 take values in open zero neighbourhoods of $\fg$ and $\fh\rtimes \fg$,
 on which the exponential functions restricts respectively to a
 diffeomorphism. Thus $\cK_{0}$ is a closed Lie subgroup of $L$.

 In the same manner, one constructs $\cK_{1}$ as a closed subgroup of
 $\cK_{0}\times C^{\infty}(\cP_{0},H\rtimes G)^{G}\times \cK_{0}$ with
 $(F,\alpha ,G)\in \cK_{1}$ if and only if $\alpha \from F\Rightarrow G$
 is a smooth natural equivalence (cf.\ Remark
 \ref{rem:equivarianceOf2-Morphisms} and \cite[Th.\ A.1]{centrExt},
 \cite[Thm.\ 1.11]{autP} for the Lie group structure on
 $C^{\infty}(\cP_{0},H\rtimes G)^{G}$). The exponential functions on
 $\cK_{0}$ and $C^{\infty}(\cP_{0},H\rtimes G)^{H\rtimes G}$ induce an
 exponential function on $\cK_{1}$ and as before, $\cK_{1}$ is a closed
 Lie subgroup. The structure maps are given by projections, embeddings
 and push-forwards by Lie group morphisms and thus they all are
 morphisms of locally convex Lie groups.
\end{proof}

\begin{corollary}
 If we endow the 2-group $\cC^{\infty}(\cP,\cG_{\Ad})^{\cG}$ with the Lie
 2-group structure from the previous theorem, then the natural action
 turns $\cP$ into  a smooth $\cC^{\infty}(\cP,\cG_{\Ad})^{\cG}$-2-space.
\end{corollary}

\begin{proof}
 This follows from the fact that evaluation maps are smooth in the
 $C^{\infty}$-topology.
\end{proof}

\begin{remark}\label{rem:crossedModuleFromGauge2-Group}
 Taking Remark
 \ref{rem:equivarianceOf2-Morphisms} into account, one obtains that
 $\cC^{\infty}(\cP,\cG_{\Ad})^{\cG}$ is associated to the smooth crossed
 module
 $(\alpha_{*},\beta_{*},C^{\infty}(\cP_{0},H)^{G}, \cMor(\cP,\cG)^{\cG})$
 with
 \begin{gather*}
  \beta_{*}(\ol{\alpha})=(\beta \circ \ol{\alpha},(\ol{\alpha}\circ
  t)\cdot (\ol{\alpha} \circ s)^{-1})\\
  \left(\alpha_{*}(\gamma_{0},\gamma_{1}).\ol{\alpha} \right)(p_{0})=
  \gamma_{0}(p_{0}).\ol{\alpha}(p_{0}).
 \end{gather*}
\end{remark}

\begin{remark}
 In
 \cite{gomi06CentralExtensionsOfGaugeTransformationGroupsOfHigherAbelianGerbes},
 there are constructed central extensions of gauge groups of (higher)
 abelian gerbes by the use of the cup-product in smooth Deligne
 (hyper-)cohomology $H^{n+2}(M,\Z^{\infty}_{D}(n+2))$. There the term
 gauge transformation is used for $H^{n+1}(M,\Z^{\infty}_{D}(n+1))$. It
 would be very interesting to explore the connection to our approach in
 order to get more general central extensions, for
 $\cC^{\infty}(\cP,\cG)^{\cG}$ in the non-abelian case (cf.\ Example
 \ref{ex:lieGroupBundle} and \cite{centrExt}).
\end{remark}

\section{Examples}

In this section, we provide some classes of examples of principal
2-bundles. This first example is an analogous construction
of the simply connected cover of a connected manifold $M$ as a 
$\pi_{1}(M)$-principal bundle. It constructs for a simply connected manifold $N$ a
principal $B \pi_{2}(N)$-2-bundle, where $B \pi_{2}(N)$ is the (discrete) 2-group
associated to the crossed module $\pi_{2}(N)\to \{*\}$. For brevity,
we shall restrict to the case where the manifold $N$ actually is a Lie group
$G$ (not necessarily finite-dimensional, so $\pi_{2}(G)\neq 0$ in
general). It already appeared implicitly at many places in the
literature (e.g.\ in
\cite{BrylinskiMcLaughlin94The-geometry-of-degree-four-characteristic%
-classes-and-of-line-bundles-on-loop-spaces.-I} and
\cite{Iglesias95La-trilogie-du-moment}) but, as far as the author knows,
it has not been worked out in terms of principal 2-bundles. It has the
correct universal property for calling it \emph{the} 2-connected cover
of $G$ (cf.\ \cite{PorstWockel08Higher-conneced-covers-of-topological%
-groups-via-categorified-central-extensions}).

\begin{example}
 Let $G$ be a 1-connected Lie group. For each $g\in G$ we choose a
 continuous path from $\gamma _{g}$ from $e$ to $g$, such that
 $\gamma_{e}\equiv e$ and $\gamma_{g}$ depends continuously on $g$ on
 some unit neighbourhood. Moreover, since $G$ is 1-connected, we find
 for each pair $(g,h)\in G^{2}$ a continuous map
 $\eta_{g,h}\from \Delta^{2}\to G$ with
 \begin{equation*}
  \partial \eta_{g,h}=\gamma_{g} + g.\gamma_{g^{-1}h}-\gamma_{h},
 \end{equation*}
 where the sum on the right is taken in the group of singular 1-chains
 of $G$. Again, we assume $\sigma_{e,e}\equiv e$ and that $\sigma_{g,h}$
 depends continuously on $g,h$ on some unit neighbourhood $U$ of $G$.
 With these choices we now set
 \begin{multline*}
  \eta _{g,h,k}\from gV \cap hV\cap kV \to \pi_{2}(G),\quad
  x\mapsto [ \underbrace{ 
  \sigma_{e,g,h}+\sigma_{e,h,k}-\sigma_{e,g,k}+
  g.\sigma_{e,g^{-1}h,g^{-1}k}}%
  _{\tx{tetrahedron with vertices }e,g,h,k}]+\\
  \underbrace{[g.(
  \sigma_{e,g^{-1}h,g^{-1}k}+
  \sigma_{e,g^{-1}k,k^{-1}x}-
  \sigma_{e,g^{-1}h,g^{-1}x}+
  g^{-1}h.\sigma_{e,h^{-1}k,h^{-1}x})}%
  _{\tx{tetrahedron with vertices }g,h,k,x}],
 \end{multline*}
 where $V\se U$ is an open symmetric unit neighbourhood with
 $V^{2}\se U$.
 If $gV \cap hV\cap kV\neq \emptyset$, then $g^{-1}h$,
 $g^{-1}k$ and $h^{-1}k$ are elements of $V^{2}\se U$. Since $\sigma$ is
 continuous on $U$ the value of $\eta_{ghk}$ does not depend on $x$, and
 $\eta _{ghk}$ is constant and in particular smooth. It is easily
 verified from the above presentation that
 $\eta_{ghk}+\eta_{gkl}=\eta_{ghl}+\eta_{hkl}$ (note that this follows
 from the above formula only because we added the second tetrahedron,
 which does not contribute to the class in $\pi_{2}(G)$). Since
 $(V_{g})_{g\in G}$ with $V_{g}:=gV$ covers $G$, $(V_{g},*,\eta_{g,h,k})$
 defines a \v{C}ech cocycle with values in $\pi_{2}(G)$ and thus a principal
 $B \pi_{2}(G)$-2-bundle, where $B\pi_{2}(G)$ denotes the 2-group
 associated to the crossed module $\pi_{2}(G)\to \{*\}$. It is fairly
 easy to check that different choices of the above data lead to
 cohomologous cocycles.
\end{example}

The remainder of this section deals with (lifting) bundle gerbes.

\begin{example}\label{ex:liftingGerbes}
 We briefly recall (abelian) bundle gerbes as introduced in
 \cite{murray96BundleGerbes} (cf.\
 \cite{aschieriCantiniLuigiBranislav05NonabelianBundleGerbesTheirDifferentialGeomoetryAndGaugeTheory},
 \cite{schweigertWaldorf07GerbesAndLieGroups}). The class of gerbes that
 connect most naturally with our principal 2-bundles are lifting gerbes, so
 we will restrict to this class (the general case can easily be adapted). Let
 $(\alpha ,\beta ,H,G)$ be a smooth crossed module and
 $\pi \from P\to M$ be a principal $G$-bundle. Moreover, we assume that
 $A:=\ker (\beta)\hookrightarrow H\twoheadrightarrow G$ is a central
 extension (i.e., we assume $\beta$ to be surjective, cf.\ \cite[Sect.\
 3]{neeb06nonAbelianExtensions}).

 There is a canonical map $f\from P\times_{M}P \to G$, determined by
 $p=p'\cdot f(p,p')$ and we consider the pull-back principal $A$-bundle
 $Q:=f^{*}(H)$ over $P\times_{M}P$.  The question that one is interested
 in is whether the $A$-action on $Q$ extends to an $H$-action, turning
 $Q$ into a principal $H$-bundle over $M$ (cf.\
 \cite{wagemannLaurent-Gengoux06},
 \cite{neeb06nonAbelianExtensionsOfTopologicalLieAlgebras}).

 With $\cG=(G,H\rtimes G)$, one can cook up a principal $\cG$-2-bundle
 $\cP$ as follows. We define objects and morphisms by
 \[
  \Obj(\cP)=P\quad \Mor(\cP)=Q,
 \]
 where we identify $Q$ with
 $\{((p,p'),h)\in P\times P\times G:p=p'\cdot \beta (h)\}$. Source and
 the target map are given by $s((p,p'),h)=p'$, $t((p,p'),h)=p$, and
 composition of morphisms is then defined by
 $((p'',p),h')\circ ((p,p'),h)=((p'',p'),h\cdot h')$ (the order of $h$
 and $h'$ is important for $((p'',p'),h\cdot h')$ to be in $Q$ again).
 In order to match this with the bundle gerbes defined in
 \cite{murray96BundleGerbes}, note that this may also be viewed as a
 $A$-bundle morphism
 \[
  \pi_{13}^{*}(Q)\to \pi_{12}^{*}(Q)\times \pi_{23}^{*}(Q)/A,
 \]
 where $A$ acts on the right hand side via the embedding
 $A\hookrightarrow A\times A$, $a\mapsto (a,a^{-1})$ and
 $\pi_{ij}\from P\times_{M}P\times_{M}P\to P\times_{M}P$ are the various
 possible projections.

 Fixing local trivialisations
 $\pi \times g_{i}\from \left.P\right|_{U_{i}}\to U_{i}\times G$ for an
 open cover $(U_{i})_{i\in I}$ defines the functors
 $\Phi_{i}\from \left.\cP\right|_{U_{i}}\to U_{i}\times \cG$ on objects
 and on morphisms we set
 $\Phi_{i}((p,p'),h)=(\pi (p),(g_{i}(p').h,g_{i}(p')))$.  The action of
 $\cG$ on $\cP$ is given by the given $G$-action on objects and by
 \[
  (((p,p'),h'),(h,g))\mapsto ((p\cdot \beta (h)\cdot g,p'\cdot
  g),g^{-1}.(h'\cdot h))
 \]
 on morphisms.

 Note that the band of this bundle is the trivial bundle over $M$,
 because $\beta$ is surjective. In particular, it has nothing to do with
 the apparent bundle $P$, which serves only as a meaningless
 intermediate space. The outer action $G/\beta (H)\to\Out(H)$ is trivial
 and lifting gerbes are classified by $H^{2}(M,A)$ (cf.\ Remark
 \ref{rem:classificationForCentralExtensions}).
\end{example}

The following example illustrates the close interplay between groups of
sections in Lie group bundles (cf.\ \cite{centrExt}) and gauge groups of
principal bundles (cf.\ \cite{autP}).

\begin{example}\label{ex:lieGroupBundle}
 Consider the crossed module $\cAut(H)$, given by the conjugation
 morphism $H\to \Aut(H)$ and the natural action of $\Aut(H)$ on $H$.
 Assuming that $H$ is finite-dimensional and $\pi_{0}(H)$ is finitely
 generated, $\Aut(H)$ becomes a Lie group, modelled on $\der(\fh)$ (cf.\
 \cite{bourbakiLie}). Moreover, we assume that we are given a smooth
 action $\lambda \from G\to \Aut(H)$ for Lie group $G$, such that the
 action lifts to a homomorphism $\varphi \from G\to H^{\flat}:=H/Z(H)$
 (this is the case, e.g., if $G$ is connected and $\fh $ is semi-simple,
 for then the derived action lifts by \cite[Prop.\
 II.6.4]{helgason78DifferentialGeometryLieGroupsAndSymmetricSpaces}).

 A Lie group bundle now arises from a principal $G$-bundle $P\to M$ by
 taking the associated bundle $P\times_{G}H$. The group of sections of
 this bundle is isomorphic to the equivariant mapping group
 $C^{\infty}(P,H)^{G}$. Considering the associated principal
 $H^{\flat}$-bundle $P^{\flat}=P\times_{\varphi}H^{\flat}$, one can ask
 whether this principal $H^{\flat}$-bundle lifts to a principal
 $H$-bundle $P^{\sharp}$. In this case, one has
 \[
  C^{\infty}(P^{\sharp},H)^{H}\cong
  C^{\infty}(P^{\flat},H)^{H^{\flat}}\cong C^{\infty}(P,H)^{G}
 \]
 (as one can see in local coordinates), and the group of sections is
 actually the gauge group of $P^{\sharp}$.

 In general, we associate to $\lambda$ the pull-back central extension
 $\varphi^{*}(H)\to G$ and thus obtain a strict 2-group $\cG$. Then we
 associate to $P^{\flat}$ the principal $\cG$-2-bundle $\cP$ of the
 lifting gerbe associated to the central extension
 $Z(H)\hookrightarrow H\twoheadrightarrow H^{\flat}$ and obtain
 $\cC^{\infty}(\cP ,\cG)^{\cG}$ as its gauge 2-group. From Remark
 \ref{rem:crossedModuleFromGauge2-Group} we see that
 $\cC^{\infty}(\cP ,\cG)^{\cG}$ is associated to a crossed module
 $(\alpha_{*},\beta_{*},H_{*},G_{*})$ with
 $H_{*}=C^{\infty}(P^{\flat},H)^{H^{\flat}}\cong C^{\infty}(P,H)^{G}$
 and
 $G_{*}\leq C^{\infty}(P^{\flat},H^{\flat})^{H^{\flat}}\nobreak \times\nobreak C^{\infty}(Q,H\nobreak \rtimes\nobreak H^{\flat})^{H\nobreak \rtimes \nobreak H^{\flat}}$.
 From the compatibility with the structure maps of $\cP$ is follows that
 \[
  (\gamma_{0},(\ol{\gamma}_{1},\gamma_{0}\circ s_{\cP}))\in
  G_{*}\Leftrightarrow \ol{\gamma}_{1}\in C^{\infty}(Q,Z(H))^{H^{\flat}}
 \]
 (note that $\gamma_{0}\circ s_{\cP}=\gamma_{0}\circ t_{\cP}$, because
 $s_{\cP}(p)=t_{\cP}(p)\cdot h$ with $h\in Z(H)$) and thus
 $G_{*}\cong C^{\infty}(P^{\flat},H^{\flat})^{H^{\flat}}\times C^{\infty}(M,Z(H))$.
 Thus $\cC^{\infty}(\cP,\cG)^{\cG}$ is in general associated to the
 crossed module
 \[
  C^{\infty}(P^{\flat},H)^{H^{\flat}}\to
  C^{\infty}(P^{\flat},H^{\flat})^{H^{\flat}}\times C^{\infty}(M,Z(H)),
  \quad \gamma\mapsto (q\circ \gamma,e_{H}).
 \]
 with the obvious point-wise action.
\end{example}

\begin{example}
 An instance of the previous example is given by considering the
 extension
 $\Gau(P)\hookrightarrow \Aut(P)\twoheadrightarrow \Diff(M)_{P}$ (cf.\
 \cite{autP}) for a finite-dimensional principal $K$-bundle $P\to M$,
 defining a crossed module by the conjugation action of $\Aut(P)$. If
 $P=M\times K$ is trivial, then $H:=\Gau(P)\cong C^{\infty}(M,K)$ and
 $\Aut(P)\cong H\rtimes \Diff(M)$ and if $K$ is compact and simple, then
 $\Aut(P)$ is an open subgroup of $\Aut(H)$ (cf.\
 \cite{guendogan2007LieAlgebrasOfSmoothSections}). Then a given action
 $G \to \Aut(H)$ lifts to $H^{\flat}$, for instance, if $G$ is connected
 and the induced action $\fg\to\cV(M)$ on the base space is trivial.
\end{example}
\appendix

\section{Appendix: Differential calculus on locally convex spaces}
\label{app:differentialCalculusOnSpacesOfMappings}

We provide some background material on spaces of mappings and their Lie
group structure in this appendix.

\begin{definition}
 \label{def:diffcalcOnLocallyConvexSpaces} Let \mbox{ $X$} and \mbox{
 $Y$} be a locally convex vector spaces and \mbox{ $U\se X$} be open.
 Then \mbox{ $f\from U\to Y$} is \emph{differentiable}  or \emph{\mbox{
 $C^{1}$}} if it is continuous, for each \mbox{ $v\in X$} the
 differential quotient
 \[
  df (x).v:=\lim_{h\to 0}\frac{f (x+hv)-f (x)}{h}
 \]
 exists and if the map \mbox{ $df\from U\times X\to Y$} is continuous.
 If \mbox{ $n>1$} we inductively define \mbox{ $f$} to be \emph{\mbox{
 $C^{n}$}} if it is \mbox{ $C^{1}$} and \mbox{ $df$} is \mbox{
 $C^{n-1}$} and to be \mbox{ $C^{\infty}$} or \emph{smooth}  if it is
 \mbox{ $C^{n}$}. We say that \mbox{ $f$} is \mbox{ $C^{\infty}$} or
 \emph{smooth} if \mbox{ $f$} is \mbox{ $C^{n}$} for all \mbox{
 $n\in\N_{0}$}. We denote the corresponding spaces of maps by \mbox{
 $C^{n}(U,Y)$} and \mbox{ $C^{\infty}(U,Y)$}.

 A (locally convex) \textit{Lie group} is a group which is a smooth
 manifold modelled on a locally convex space such that the group
 operations are smooth.
\end{definition}

\begin{proposition}
 \label{prop:localDescriptionsOfLieGroups} Let \mbox{\mbox{ $G$}} be a
 group with a locally convex manifold structure on some subset
 \mbox{\mbox{ $U\se G$}} with \mbox{\mbox{ $e\in U$}}. Furthermore,
 assume that there exists \mbox{\mbox{ $V\se U$}} open such that
 \mbox{\mbox{ $e\in V$}}, \mbox{\mbox{ $VV\se U$}}, \mbox{\mbox{
 $V=V^{-1}$}} and
 \begin{itemize}
  \item [i)] \mbox{\mbox{ $V\times V\to U$}}, \mbox{\mbox{
        $(g,h)\mapsto gh$}} is smooth,
  \item [ii)] \mbox{\mbox{ $V\to V$}}, \mbox{\mbox{ $g\mapsto g^{-1}$}}
        is smooth,
  \item [iii)] for all \mbox{\mbox{ $g\in G$}}, there exists an open
        unit neighbourhood \mbox{\mbox{ $W\se U$}} such that
        \mbox{\mbox{ $g^{-1}Wg\se U$}} and the map \mbox{\mbox{
        $W\to U$}}, \mbox{\mbox{ $h\mapsto g^{-1}hg$}} is smooth.
 \end{itemize}
 Then there exists a unique locally convex manifold structure on
 \mbox{\mbox{ $G$}} which turns \mbox{\mbox{ $G$}} into a Lie group,
 such that \mbox{\mbox{ $V$}} is an open submanifold of \mbox{\mbox{
 $G$}}.
\end{proposition}

\begin{definition}\label{def:exponentialfunction}
 Let \mbox{\mbox{ $G$}} be a locally convex Lie group. The group
 \mbox{\mbox{ $G$}} is said to have an \emph{exponential function} if
 for each \mbox{\mbox{ $x \in \fg$}} the initial value problem
 \[
  \gamma (0)=e,\quad \gamma '(t)=T\lambda_{\gamma (t)}(e).x
 \]
 has a solution \mbox{\mbox{ $\gamma_{x}\in C^{\infty} (\R,G)$}} and the
 function
 \[
  \exp_{G}:\fg\to G,\;\;x \mapsto \gamma_x (1)
 \]
 is smooth. Furthermore, if there exists a zero neighbourhood
 \mbox{\mbox{ $W\se \fg$}} such that \mbox{\mbox{
 $\left.\exp_{G}\right|_{W}$}} is a diffeomorphism onto some open unit
 neighbourhood of \mbox{\mbox{ $G$}}, then \mbox{\mbox{ $G$}} is said to
 be \emph{locally exponential}.
\end{definition}

\begin{lemma}
 \label{lem:interchangeOfActionsOnGroupAndAlgebra} If \mbox{\mbox{ $G$}}
 and \mbox{\mbox{ $G'$}} are locally convex Lie groups with exponential
 function, then for each morphism \mbox{\mbox{ $\alpha :G\to G'$}} of
 Lie groups and the induced morphism \mbox{\mbox{
 $d\alpha (e):\fg \to \fg'$}} of Lie algebras, the diagram
 \[
  \begin{CD}
   G @>\alpha >>G'\\
   @AA\exp_{G}A  @AA\exp_{G'}A\\
   \fg@>d\alpha (e)>> \fg'
  \end{CD}
 \]
 commutes.
\end{lemma}

\begin{remark}
 \label{rem:banachLieGroupsAreLocallyExponential} The Fundamental
 Theorem of Calculus for locally convex spaces (cf.\ \cite[Th.\
 1.5]{gloeckner02b}) yields that a locally convex Lie group \mbox{\mbox{
 $G$}} can have at most one exponential function (cf.\ \cite[Lem.\
 II.3.5]{neeb06}).

 Typical examples of locally exponential Lie groups are Banach-Lie
 groups (by the existence of solutions of differential equations and the
 inverse mapping theorem, cf.\ \cite{lang99}) and groups of smooth and
 continuous mappings from compact manifolds into locally exponential
 groups (\cite[Sect.\ 3.2]{gloeckner02a}, \cite{smoothExt}). However,
 diffeomorphism groups of compact manifolds are never locally
 exponential (cf.\ \cite[Ex.\ II.5.13]{neeb06}) and direct limit Lie
 groups not always (cf.\ \cite[Rem.\ 4.7]{gloeckner05}). For a detailed
 treatment of locally exponential Lie groups and their structure theory
 we refer to \cite[Sect.\ IV]{neeb06}.
\end{remark}

\begin{remark}
 The most interesting examples of infinite-dimensional Lie groups for
 this article shall be groups of smooth mappings $C^{\infty}(M,G)$ from
 a compact manifold $M$ (possibly with boundary) to an arbitrary Lie
 group $G$. These groups possess natural Lie group structures if one
 endows them with the initial topology with respect to the embedding
 \[
  C^{\infty}(M,G)\hookrightarrow \prod_{k\in
  \N_{0}}C(T^{k}M,T^{k}G)_{c}, \quad \gamma\mapsto (T^{k}\gamma)_{k\in
  \N_{0}}.
 \]
 The Lie algebra is $C^{\infty}(M,\fg)$ (with the above topology), where
 $\fg$ is the Lie algebra of $G$ (all spaces are endowed with point-wise
 operations). Details can be found in \cite[Sect.\ 3.2]{gloeckner02a}
 and \cite[Sect.\ 4]{smoothExt}.
\end{remark}

\section*{Acknowledgements}

The author would like to thank Thomas Schick, Peter Arndt, Ulrich
Pennig, Sven Porst, and Alessandro Fermi for various useful discussions
and Alessandro Fermi for proof-reading the manuscript. Moreover, he
would like to thank Ieke Moerdijk for writing the excellent introduction
\cite{moerdijk02IntroductionToTheLanguageOfStacksAndGerbes}, which
clarified the subject essentially and supported many ideas of this
paper. Urs Schreiber gave a couple of very useful comments on various
formation stages, providing the background and global picture and
proof-read the paper thoroughly. Eventually, Chris Schommer-Pries
pointed out Remark \ref{rem:morphismsAndOpenCovers} and hereby helped
correcting an error in a previous version of the paper.

The work on this article was financially supported by the Graduiertenkolleg
``Gruppen und Geometrie''.

\def\polhk#1{\setbox0=\hbox{#1}{\ooalign{\hidewidth
  \lower1.5ex\hbox{`}\hidewidth\crcr\unhbox0}}} \def\cprime{$'$}

\end{document}